\newif\ifspringer
\newif\ifelsevier
\elseviertrue 

\ifspringer
\RequirePackage{fix-cm} 
\documentclass[smallcondensed,nospthms]{svjour3}
\usepackage[paperwidth=15.51cm,paperheight=23.51cm,left=1.75cm,right=1.75cm,top=2cm,bottom=2cm]{geometry} 
\fi

\ifelsevier
\documentclass[preprint,3p,times]{elsarticle}
\fi

\usepackage{widetext}
\usepackage[english]{babel}
\usepackage[utf8]{inputenc}
\usepackage{natbib}
\usepackage{pdflscape} 
\usepackage{amssymb}
\usepackage{amsmath}
\numberwithin{equation}{section} 
\usepackage{amsthm}
\usepackage{mathtools}
\ifspringer
\usepackage{bm} 
\fi
\usepackage[normalem]{ulem}
\usepackage{enumitem}
\usepackage{marvosym}
\usepackage[dvipsnames,svgnames]{xcolor}
\usepackage{array}
\usepackage{multirow}
\ifspringer
\usepackage{mathptmx} 
\fi
\ifelsevier
\usepackage[nodots]{numcompress} 
\fi

\usepackage{graphicx}
\ifspringer
\usepackage[FIGBOTCAP,TABTOPCAP,scriptsize]{subfigure} 
\fi
\ifelsevier
\usepackage[FIGBOTCAP,TABTOPCAP]{subfigure} 
\fi
\usepackage{tikz}
\usetikzlibrary{calc,shapes,arrows,backgrounds}
\tikzstyle{every picture}+=[remember picture]
\usepackage{pgfplots} 

\usepackage{hyperref}
\hypersetup{colorlinks=true,linkcolor=blue,citecolor=blue,filecolor=blue,urlcolor=blue}
\usepackage[all]{hypcap}

\mathtoolsset{showonlyrefs} 
\allowdisplaybreaks 
\ifspringer
\smartqed 
\fi
\ifelsevier
\biboptions{sort&compress,numbers}
\fi

\usepackage{xspace}
\makeatletter
\DeclareRobustCommand\onedot{\futurelet\@let@token\@onedot}
\newcommand{\@onedot}{\ifx\@let@token.\else.\null\fi\xspace}
 
\newcommand{\ie}{{i.e}\onedot}

\newcommand{\wrt}{{w.r.t}\onedot}

\makeatother

\newtheorem{rem}{Remark}

\theoremstyle{definition}
\newtheorem{defn}{Definition}
\newtheorem{exmp}{Example}
\newtheorem{alg}{Algorithm}

\newcommand{\RR}{\mathbb{R}}

\newcommand{\vectbi}[1]{\boldsymbol{#1}} 
\newcommand{\matri}[1]{#1}               
\newcommand{\vect}[1]{\vectbi{#1}}
\newcommand{\matr}[1]{\matri{#1}}

\DeclareMathOperator{\Span}{span}
\renewcommand{\leq}{\leqslant}
\renewcommand{\geq}{\geqslant}

\newcommand{\rwh}{\hat}
\newcommand{\PP}{\mathcal{P}}
\ifspringer
\DeclareMathAlphabet{\mathcalb}{OMS}{cmsy}{b}{n} 
\DeclareMathAlphabet{\mathcal}{OMS}{cmsy}{m}{n} 
\newcommand{\TTS}{\mathcalb{U}}
\newcommand{\MMS}{\mathcal{M}}
\fi
\ifelsevier
\newcommand{\TTS}{\vect{\mathcal{U}}}
\newcommand{\MMS}{\mathcal{M}}
\fi
\newcommand{\TT}{\mathcal{U}}

\newcommand{\bM}{\vect{M}}

\newcommand{\bDelta}{\vect{\Delta}}

\newlength{\casesvsep}
\newcommand{\figpath}{figure}
\graphicspath{{\figpath/}}

\makeatletter
\newcommand*{\shifttext}[2]{%
\settowidth{\@tempdima}{#2}%
\makebox[\@tempdima]{\hspace*{#1}#2}%
}
\makeatother

\newcommand{\qtext}[1]{``#1''}

\ifspringer

\fi
\ifelsevier

\fi

\ifspringer
\journalname{\dots}
\fi

\makeatletter
\providecommand{\doi}[1]{%
  \begingroup
    \let\bibinfo\@secondoftwo
    \urlstyle{rm}%
    \href{http://dx.doi.org/#1}{%
      doi:\discretionary{}{}{}%
      \nolinkurl{#1}%
    }%
  \endgroup
}
\makeatother

\begin{document}
\newcommand{\titletext}{
A practical method for computing with piecewise Chebyshevian splines
}
\newcommand{\titlerunningtext}{
A practical method for computing with piecewise Chebyshevian splines
}

\newcommand{\abstracttext}{
A piecewise Chebyshevian spline space is \emph{good for design} when it possesses a B-spline basis and this property is preserved under knot insertion.
The interest in such kind of spaces is justified by the fact that, similarly as for polynomial splines, the related parametric curves exhibit the desired properties of convex hull inclusion, variation diminution and intuitive relation between the curve shape and the location of the control points.
For a good-for-design space, in this paper we construct a set of functions, called transition functions, which allow for efficient computation of the B-spline basis, even in the case of nonuniform and multiple knots.
Moreover, we show how the spline coefficients of the representations associated with a refined knot partition and with a raised order can conveniently be expressed by means of transition functions.
This result allows us to provide effective procedures that generalize the classical knot insertion and degree raising algorithms for polynomial splines.
We further discuss how the approach can straightforwardly be generalized to deal with geometrically continuous piecewise Chebyshevian splines as well as with splines having section spaces of different dimensions.
From a numerical point of view, we show that the proposed evaluation method is easier to implement and has higher accuracy than other existing algorithms.
}

\ifspringer
\newcommand{\separ}{\and}
\fi
\ifelsevier
\newcommand{\separ}{\sep}
\fi

\newcommand{\keytext}{
(Piecewise) Chebyshevian splines \separ B-spline basis \separ Knot insertion \separ Order elevation \separ Computational algorithms \separ Transition functions.
}


\ifspringer
\title{\titletext
}
\titlerunning{\titlerunningtext} 

\author{
Carolina Vittoria Beccari \and
Giulio Casciola           \and
Lucia Romani
}
\authorrunning{C.V.~Beccari, G.~Casciola, L.~Romani} 

\institute{
C.V.~Beccari (\Letter) \at
 Department of Mathematics, University of Bologna,\\
 Piazza di Porta San Donato 5, 40126 Bologna, Italy\\
 \email{carolina.beccari2@unibo.it} 
\and
G.~Casciola \at
 Department of Mathematics, University of Bologna,\\
 Piazza di Porta San Donato 5, 40126 Bologna, Italy\\
 \email{giulio.casciola@unibo.it} 
\and
L.~Romani \at
 Department of Mathematics, University of Bologna,\\
 Piazza di Porta San Donato 5, 40126 Bologna, Italy\\
 \email{lucia.romani@unibo.it} 
}

\date{Received: date / Accepted: date}

\maketitle

\begin{abstract}
\abstracttext
\end{abstract}
\fi

\ifelsevier
\begin{frontmatter}

\title{\titletext}

\author[label2]{Carolina Vittoria Beccari\corref{cor1}}
\ead{carolina.beccari2@unibo.it}
\author[label2]{Giulio Casciola}
\ead{giulio.casciola@unibo.it}
\author[label2]{Lucia Romani}
\ead{lucia.romani@unibo.it}

\cortext[cor1]{Corresponding author.}

\address[label2]{Department of Mathematics, University of Bologna,
Piazza di Porta San Donato 5, 40126 Bologna, Italy}

\begin{abstract}
\abstracttext
\end{abstract}


\end{frontmatter}
\fi

\section{Introduction}
\label{sec:intro}
Extended Chebyshev (EC) spaces represent a natural generalization of polynomials.
They may contain transcendental functions (and thus are capable of reproducing circles, ellipses and many other curves that cannot be represented by polynomials) and in addition they provide parameters which can be
exploited to control the behavior of parametric curves and to accomplish shape-preserving
approximations.
Functions having pieces in different EC-spaces are called \emph{piecewise Chebyshevian splines}.
Classical polynomial splines, Chebyshevian splines having all pieces in the same EC-space (such as, e.g., hyperbolic splines) and generalized splines
\cite{Cos2000}
or GB-splines \cite{KS1999} are all special instances in this wide class.
Piecewise Chebyshevian splines are of great interest for their ability to combine the local nature of splines with the
diversity of shape effects provided by the wide range of known EC-spaces.
Precisely, when the underlying section spaces are suitably chosen, they allow for the exact representation of all possible varieties of transcendental curves. For example, segments of helices, which can be used for
tool path description, simulation of kinematic motion and design of highways, can be represented accurately with a combination of trigonometric functions and polynomials, whereas they do not admit an exact representation in terms of polynomials
or rational polynomials.
The interest in piecewise Chebyshevian splines is also demonstrated by their numerous applications, ranging from Approximation Theory and Geometric Design \cite{Bracco2014a,Gonz2016,RSA2014,RM2019}, to Subdivision Theory and Multiresolution Analysis
\cite{CGR2016,Jeong2013647,Romani2021,KU2005,LY2011,LM2012},
Image Processing \cite{BSUU16}, numerical solution of PDEs \cite{YoonWENO2016,Yuan2006295,Zhu2013} and Isogeometric Analysis \cite[pp.1-69]{BS2012}, \cite{MPS2011b,Manni_iga1}.

In the following, for the sake of shortness, we will omit the word \qtext{piecewise} and just refer to these splines as Chebyshevian splines, as is the case for several other previous papers (see, e.g., \cite{SPE2020}).

\subsection{State-of-the-art and new contributions on piecewise Chebyshevian splines (for short, Chebyshevian splines)}

Chebyshevian splines, considered so far in the literature, can be either parametrically or geometrically continuous. The latter are piecewise functions where the continuity conditions between adjacent segments are expressed by means of connection matrices \cite{Barry1996}.
In both cases, to be of interest for Geometric Design or other applications,
a Chebyshevian spline space should not only possess a B-spline basis
-- in the usual sense of a normalised basis composed of minimally supported splines --
but this should be true also for any other space obtained by knot insertion.
This requirement is equivalent to existence of blossoms in the space, see \cite{Maz2011b} and references therein, and, when it holds, it is said that the space is \emph{good for design}, or has a \emph{refinable} B-spline basis.
As a consequence of blossoms, the B-spline basis is the (unique) optimal normalized totally positive basis and all the classical geometric design algorithms can be developed in the space.

Theoretical results to state whether or not a specific (and usually simple, e.g., with low dimension and maximum continuity) space is good for design have been provided (see \cite{Maz2011b,Mazure2018550}
and references therein), as well as numerical procedures which work for every arbitrary space \cite{BCM2016,BCM2019}.

Regarding the construction and evaluation of Chebyshevian splines, there exist a very general method \cite[Section 7]{Maz2011b} that, even if can be applied to any arbitrary Chebyshevian spline space, is unfortunately of limited practical use since, on one side, it relies on the so-called \emph{weight functions} of the spline space (which are generally not known) and, on the other side, its integral and recursive nature makes it unsuited for numerical computations. Thus, the problem of efficiently and accurately evaluating Chebyshevian splines has recently returned to the fore. In particular, the paper \cite{AGO2019} has proposed a method for the evaluation of Chebyshevian Bernstein bases, while the two papers \cite{chebMD2019,SPE2020} have focussed on the evaluation of splines (piecewise functions), also dealing with the extremely general case in which section spaces may have different dimensions (which they call \emph{multi-degree Tchebycheffian splines}).

In the wake of such recent works, our paper wants to propose an alternative approach for the construction and efficient evaluation of Chebyshevian B-spline bases, which produces accurate results, being at the same time simple from a conceptual point of view and easy to implement.
Our algorithm exploits the fact that, if a space is good for design, then any Hermite interpolation problem (under proper conditions relating nodes and knots) has a unique solution. As a result, the algorithm involves the resolution of a certain number of small-size linear systems (as many as the spline space dimension minus one). These linear systems do not lead directly to the calculation of the B-spline basis, but rather generate a special set of functions that are called \emph{transition functions}. The latter form a basis of the space of interest such that each B-spline function is the difference of two successive transition functions.

As a historical note, it is interesting to remark that the idea of transition functions, in its first
form, had already appeared in the original work by Pierre B\'ezier as a method for evaluating the Bernstein basis.
Much more recently, transition functions were used in connection with the generation of local piecewise polynomial interpolants \cite{ABC2013a} and also in the construction of the B-spline basis for \emph{multi-degree} splines \cite{BCM2017} (splines whose pieces are polynomials of different degrees), as well as for the analysis of the structural properties of such spaces \cite{TSHMH2018}.

The possibility of constructing transition functions for Chebyshevian splines was originally mentioned in our technical report \cite{BCR2019_comp} that was registered in the arXiv repository even before the publication of the above cited papers.
At that time, what had held us back from finalizing the idea of using the transition functions for the evaluation of Chebyshevian splines was the fact that sometimes the computation of the functions that generate these spaces (think, e.g., of hyperbolic functions) can be subject to numerical problems that inevitably affect the accuracy of any algorithm for the evaluation of their B-spline basis. The most recent works on the subject highlight that these well-known problems have not actually hindered the interest and diffusion of Chebyshevian splines and have therefore convinced us that evaluation by means of transition functions is a viable path and has strong practical interest. In particular, the only alternative evaluation method \cite{chebMD2019,SPE2020} requires more elaborate calculations and is subject to similar, and in some cases even worse, numerical inaccuracy issues, as we will show in the last section of this paper.
Our algorithm is general (it works for any Chebyshevian spline space), conceptually simple, and easy to implement. It has the advantage of being local (namely it involves calculations with small-size variables and a localized change in the space does not involve recalculating the basis in full, but only the modified functions) and is also accurate in all situations of common interest. Our method also manages the case in which the continuity conditions are given by connection matrices (in the literature referred to as \emph{geometric continuity}), which is not included in \cite{chebMD2019,SPE2020}. In describing our algorithm we will focus on splines where the section spaces are all of the same dimension. This choice, in addition to simplifying the presentation, is natural as only in this case there are theoretical and computational methods to ensure that a space is good for design, or, in other words, to ensure the existence of the refinable B-spline basis that we build.
The generalization of our computational proposal to spaces with sections of different dimensions is straightforward, and all the necessary details towards this generalization will be provided at the end of this paper together with a comparison with the algorithm in \cite{chebMD2019,SPE2020}.

More precisely, the remainder of the paper is organized as follows.
In Section \ref{sec:blending} we introduce the transition functions and discuss how to use them for the evaluation of the B-spline basis of a Chebyshevian spline space. In Section \ref{sec:knotins_ordelev} we exploit these functions to develop the classical tools of knot insertion and dimension elevation for Chebyshevian splines.
Section \ref{sec:gc_extension} discusses how the transition functions can be used in the more general setting of geometrically continuous splines. The generalization of our computational proposal to splines with sections of different dimensions is presented in Section \ref{sec:mo-spline}. Finally Section \ref{sec:comparison} presents a comparison with the algorithm in \cite{chebMD2019,SPE2020}, which is the only equally general evaluation method available at present.

\section{Transition functions and B-spline bases for Chebyshevian spline spaces}
\label{sec:blending}

In this section we introduce the definition of transition functions and exploit such functions to efficiently compute the B-spline basis of any Chebyshevian spline space good for design, regardless of the kind and dimension of the different section spaces.

\subsection{Preliminary notions on EC-spaces, Chebyshevian splines and B-spline bases}
\label{sec:preliminaries}
Let us proceed by recalling the definition of EC-spaces, which are the building blocks of the splines we are interested in.

\begin{defn}[Extended Chebyshev space]\label{defn:ECspace}
An $m$-dimensional space $\TT$ contained in $C^{m-1}(I)$, $I\subset \RR$, is an \emph{Extended Chebyshev space} (for short, EC-space) on $I$ if any nonzero element of $\TT$ vanishes at most $m-1$ times in $I$, counting multiple zeros as far as possible for $C^{m-1}$ functions (that is, up to $m$), or, equivalently, if any Hermite interpolation problem in $m$ data in $I$ has a unique solution in $\TT$.
\end{defn}

Subdividing the domain interval into a number of consecutive subintervals, it is possible to introduce the definition of Chebyshevian splines, referred to as piecewise functions whose pieces may belong to different EC-spaces.
For the moment, we limit ourselves to considering Chebyshevian splines where adjacent spline pieces are connected via the standard parametric continuity and postpone the rather straightforward generalization of the proposed evaluation method to geometrically continuous splines to a later section (specifically Section \ref{sec:gc_extension}).

To define a Chebyshevian spline space we need to introduce the following quantities.
Let $[a,b]$ be a bounded and closed interval, and $\bDelta\coloneqq\left\{x_i\right\}_{i=1,\dots,q}$ be a set of points such that $a \equiv x_{0} < x_{1} < \ldots <x_{q} < x_{q+1} \equiv b$, which determine a partition of $[a,b]$ into subintervals $\left[x_i,x_{i+1}\right]$, $i=0,\dots,q$.
Moreover, let $m$ be a positive integer, and $\MMS\coloneqq\left(\mu_1,\dots,\mu_q\right)$ be a vector of positive integers, which we call \emph{multiplicities}, such that $0 \leq \mu_i < m$ for every $i=1,\dots,q$.
We denote by $\TTS_m\coloneqq\{\TT_{0,m},\dots,\TT_{q,m}\}$ a sequence of spaces of dimension $m$ such that $\TT_{i,m}$ is an EC-space containing constants on $[x_i,x_{i+1}]$ and such that $D\TT_{i,m}\coloneqq\{Du \coloneqq u^{\prime} \; | \; u \in \TT_{i,m}\}$ is an $(m-1)$-dimensional EC-space on $[x_i,x_{i+1}]$, $i=0,\dots,q$. These requirements guarantee that $\TT_{i,m}$ has the Bernstein basis both on the interval $[x_i,x_{i+1}]$ and on all subintervals contained in it. The supremum of all real numbers $L$ such there exist a Bernstein basis on and only on compact intervals of length less than $L$ is called the \emph{critical length for design} and can be computed by either theoretical tools or numerical methods (see \cite{BCM2020a} and references therein).

\begin{defn}[Piecewise Chebyshevian splines (for short, Chebyshevian splines)]\label{def:QCS}
We define the set of \emph{piecewise Chebyshevian splines} (for short, Chebyshevian splines) of order $m$ with break-points $x_1,\dots,x_q$ of multiplicities $\mu_1,\dots,\mu_q$ as
\begin{align}
S(\TTS_m,\MMS,\bDelta) \coloneqq & \left\{ s \,\big|\, \mbox{there exist } s_i\in\TT_{i,m}, i=0,\dots,q, \right.\\[-1ex]
& \hspace{-1.8cm} \mbox{ such that:} \\[-1ex]
& \hspace{-2cm} \begin{minipage}[b]{0.95\linewidth}
\begin{enumerate}[label=\roman*)]
\item $s(x)=s_i(x)$ for $x\in [x_i,x_{i+1}], i=0,\dots,q$;
\item $D^\ell s_{i-1}(x_i)=D^\ell s_{i}(x_i)$ for $\ell=0,\dots,m-\mu_i-1, \quad i=1,\dots,q \left. \vphantom{\big|} \right\}.$
\end{enumerate}
\end{minipage}
\end{align}
\end{defn}
In consequence of standard arguments $S(\TTS_m,\MMS,\bDelta)$ is a function space of dimension $m+K$ with $K\coloneqq\sum_{i=1}^q \mu_i$.

\begin{defn}[Extended partition]\label{def:part_estesa}
The set of knots $\bDelta^*\coloneqq\left\{t_i\right\}_{i=1,\dots,2m+K}$, with $K=\sum_{i=1}^q \mu_{i}$, is called an \emph{extended partition} associated with $S(\TTS_m,\MMS,\bDelta)$ if and only if:
\begin{enumerate}[label=\roman*)]
 \item $t_1 \leq t_2 \leq \dots \leq t_{2m+K}$;
 \item $t_m \equiv a$ and $t_{m+K+1} \equiv b$;
 \item $\left\{t_{m+1}, .., t_{m+K}\right\} \equiv
        \{ \underbrace{x_1,..,x_1}_{\mu_1 \text{ times}}, \dots, \underbrace{x_q,..,x_q}_{\mu_q \text{ times}} \}$,
\end{enumerate}
{ with $t_i$, $i<m$ and $i>m+K+1$ arbitrarily chosen so as to satisfy i).}
\end{defn}

In the following we use the terms \emph{break-point}, \emph{knot} and \emph{node} in the following sense: a break-point
is the junction between two spline segments, a knot corresponds to a break-point repeated as many times as its multiplicity, and a node is the abscissa of an interpolation point. We also introduce the notation
\begin{equation}\label{eq:muLR}
\begin{array}{l}
\mu_i^\mathrm{L}\coloneqq\max \{ j \geq 0 \ | \ t_{i-j}=t_i \}+1,\\[1ex]
\mu_i^\mathrm{R}\coloneqq\max \{ j\geq 0 \ | \ t_i=t_{i+j} \}+1,
\end{array}
\end{equation}
to identify the left ($\mu_i^\mathrm{L}$) and right ($\mu_i^\mathrm{R}$) multiplicity associated with a knot $t_i$.

For simplicity and without loss of generality, we will confine our discussion to the case of extended partitions
without external additional knots, \ie assuming $x_0 \equiv a$, $x_{q+1} \equiv b$ and $\mu_0=\mu_{q+1}=m$.
Chebyshevian splines can similarly be constructed from a general partition with external knots and break-points. In the latter case, we shall assign an auxiliary EC-space to each interval which is not contained in $[a,b]$, in such a way that the continuity conditions between adjacent spline pieces are well-defined.

\noindent
The next definition characterizes the B-spline basis of a Chebyshevian spline space in terms of its distinguishing properties (see also \cite{Maz2011b}). { We equivalently refer to it as either the \emph{Chebyshevian B-spline basis} or simply the B-spline basis, since it has the same properties as the polynomial B-spline basis.}

\begin{defn}[B-spline basis]\label{def:B-spline basis}
A sequence $\{N_{i,m}\}_{i=1,\dots,m+K}$ of elements of $S(\TTS_m,\MMS,\bDelta)$ is the \emph{B-spline basis} if it meets the following requirements:
\begin{enumerate}[label=\roman*)]
\setlength{\itemsep}{4pt}
\item \label{propty:BS1} \emph{support property}: for each $i \in \{1,\dots,m+K\}$, \ $N_{i,m}(x)=0$ for $x \notin [t_i, t_{i+m}]$;
\item \label{propty:BS2} \emph{positivity property}: for each $i \in \{1,\dots,m+K\}$, \ $N_{i,m}(x) > 0$ for $x \in (t_i,t_{i+m})$;
\item \label{propty:BS3} \emph{endpoint property}: for each $i \in \{1,\dots,m+K\}$, \ $N_{i,m}$ vanishes exactly $m-\mu_i^\mathrm{R}$ times at $t_i$ and exactly $m-\mu_{i+m}^\mathrm{L}$ times at $t_{i+m}$ where, according to \eqref{eq:muLR},
$\mu_i^\mathrm{R}\coloneqq\max \{ j\geq 0 \ | \ t_i=t_{i+j} \}+1$ and $\mu_{i+m}^\mathrm{L}\coloneqq\max \{ j \geq 0 \ | \ t_{i+m-j}=t_{i+m} \}+1$;
\item \label{propty:BS4} \emph{normalization property}:
$$\displaystyle \sum_i N_{i,m}(x) = 1, \; \forall x \in [a,b].$$
\end{enumerate}
\end{defn}

Any spline $s \in S(\TTS_m,\MMS,\bDelta)$ has an expansion in the B-spline basis $N_{i,m}$, $i=1,\dots,m+K$, of the form
\begin{equation}\label{eq:spline_global}
s(x)=\sum_{i=1}^{m+K} c_i\,N_{i,m}(x), \qquad x\in[a,b].
\end{equation}
Furthermore, in virtue of the local support of the basis functions, there are only $m$ nonvanishing B-spline functions on a single interval,  and the above relation simplifies to
\begin{equation}\label{eq:spline}
s(x)=\sum_{i=\ell-m+1}^{\ell} c_i\,N_{i,m}(x), \qquad x\in[t_\ell,t_{\ell+1}).
\end{equation}

As already stated in the introduction, to be of interest for Geometric Design or other applications, a Chebyshevian
spline space should not only possess the B-spline basis, but it should also be \emph{good for design} (see \cite{Maz2011b}) in the sense specified in the next definition.

\begin{defn}[Good for design spline space]\label{def:GFD}
A Chebyshevian spline space is \emph{good for design} if it possesses the B-spline basis and so does any spline space obtained from it by knot insertion.
\end{defn}

In view of the above definition, we will also say that a good for design spline space has a refinable B-spline basis.  As proved in \cite{Maz2011b} (see also references therein) a space is good for design if and only if  blossoms exist for that space. The existence of blossoms guarantees that a Chebyshevian spline space \qtext{behaves like} a polynomial space and, in particular, that parametric curves in this space exhibit the desired properties of convex hull inclusion, variation diminution and the intuitive relation between the curve shape and the location of the control points (for further acquaintance see, e.g., \cite{Maz2011b,LM2012} and references therein).
Moreover, in a good for design spline space any Hermite interpolation problem has a unique solution under suitable interlacing conditions \cite{BMuhl2003}.

In the literature, several examples of good for design Chebyshevian spline spaces can be found,  and these have been dubbed under different names including Unified Extended splines \cite{WF2008},
NUAHT-splines \cite{XW2007}, cycloidal splines \cite{BR07}, exponential splines \cite{Unser2005} and generalized splines \cite{BMuhl2003,CMP2003,MP2010}. For many of these spaces, the property of being good for design comes directly
from the fact that their B-spline basis is constructed by means of an integral recurrence relation.
Furthermore, thanks to the computational methods recently proposed in \cite{BCM2016,BCM2019}, it is possible to establish if any arbitrary Chebyshevian spline space is good for design by efficient and accurate numerical procedures.

\subsection{Construction of the B-spline basis of Chebyshevian spline spaces via transition functions}
\label{sec:construction}
In this section we introduce the transition functions for the space $S(\TTS_m,\MMS,\bDelta)$ and show how to use these functions to construct and evaluate the B-spline basis.

\begin{defn}[Transition functions]\label{def:trans_function}
Let $S(\TTS_m,\MMS,\bDelta)$ be a Chebyshevian spline space of dimension $m+K$ with $\bDelta=\{x_i\}_{i=1,\dots,q}$ a partition of $[a,b]$ and $\bDelta^*=\{t_i\}_{i=1, \dots, 2m+K}$ the associated extended partition.
Let also suppose that $S(\TTS_m,\MMS,\bDelta)$ has a B-spline basis
$\{N_{i,m}\}_{i=1,\dots,m+K}$.
We call \emph{transition functions} the piecewise functions $f_{i}$, given by:
\begin{equation}\label{eq:f_i_N_i}
 f_{i}=\sum_{j=i}^{m+K} N_{j,m}, \quad i=1,\dots,m+K.
\end{equation}
\end{defn}

\noindent
One can immediately verify that
the transition functions are a basis for $S(\TTS_m,\MMS,\bDelta)$ by observing that relations \eqref{eq:f_i_N_i} can be written in the matrix form
\begin{equation}\label{eq:f_i_N_i_matrix}
\left( \begin{array}{c}
f_{1} \\ f_{2} \\ \vdots \\ \vdots \\ f_{m+K} \\
\end{array} \right) =
\left( \begin{array}{cccc}
1 & 1 & \cdots & 1\\
0 & 1 & \ddots & \vdots \\
\vdots & \ddots & \ddots  & 1 \\
0 & \cdots & 0 & 1 \\
\end{array} \right)
\left( \begin{array}{c}
N_{1,m} \\ N_{2,m} \\ \vdots \\ \vdots \\ N_{m+K,m} \\
\end{array} \right)
\end{equation}
with a nonsingular matrix.
Moreover, assuming $f_{m+K+1}\equiv 0$ and inverting relation \eqref{eq:f_i_N_i_matrix}, we can express the B-spline basis in terms of transition functions as
\begin{equation}\label{eq:N}
N_{i,m} = f_{i}-f_{i+1}, \quad i=1,\dots,m+K.
\end{equation}

From the properties of the B-spline basis listed in Definition \ref{def:B-spline basis} we can deduce that $f_1(x) = 1$, for all $x \in [a,b]$, and that the piecewise functions $f_i$, $i=2,\dots,m+K$, have the following characteristics:
\begin{enumerate}
\item[TF-A.] \label{propty:f1}
Each transition function is nontrivial in a closed and bounded interval and precisely
\begin{equation}
f_{i}(x)=
\begin{cases}
0, & x \leq t_i,\\
1, & x \geq t_{i+m-1};
\end{cases}
\end{equation}

\item[TF-B.] \label{propty:f2} If $t_{i}<t_{i+m-1}$, then
$f_{i}$ vanishes $m-\mu_i^\mathrm{R}$ times at $t_i$ and $1-f_{i}$ vanishes $m-\mu_{i+m-1}^\mathrm{L}$ times at $t_{i+m-1}$.
In particuar, the following bounds on the number of zero derivatives hold
\begin{equation}
\begin{array}{l}
D_+^{r} f_{i}(t_{i}) = 0, \quad r=0,\dots,k_{i}^R, \\
\\
D_+^{k_{i}^R+1}f_{i}(t_{i}) > 0,\\
\end{array}
\qquad \quad
\hbox{and}
\qquad \quad
\begin{array}{l}
D_-^{r} f_{i}(t_{i+m-1}) = \delta_{r,0}, \quad r=0,\dots,k_{i+m-1}^L, \\
\\
(-1)^{k_{i+m-1}^L}D_-^{k_{i+m-1}^L+1}f_{i}(t_{i+m-1}) > 0, \\
\end{array}
\end{equation}
where
\begin{equation}\label{eq:ki}
k_{i}^R \coloneqq m-\mu_i^\mathrm{R}-1 \quad \textrm{and} \quad k_{i+m-1}^L \coloneqq m-\mu_{i+m-1}^L-1.
\end{equation}
\end{enumerate}

\begin{rem}[Bernstein basis]\label{rem:Bernst}
A Chebyshevian spline space with empty knot partition (\ie, $\bDelta=\emptyset$) is just an EC-space and formula \eqref{eq:N} yields its Bernstein basis.
In this case, we denote the basis functions of the $m$-dimensional EC-space by $B_{0,m-1},\dots,B_{m-1,m-1}$, and the corresponding transition functions by $g_0,\dots,g_{m-1}$. In particular, formula \eqref{eq:N} reads now as
\begin{equation}\label{eq:B}
B_{i,m-1} = g_i-g_{i+1}, \quad i = 0,\dots,m-1,
\end{equation}
where $g_0\equiv 1$, $g_{m}\equiv 0$ and the functions $g_{i}$, $i=1,\dots,m-1$ are determined by the conditions $D^r g_{i}(a)=0$, \, $r=0,\dots,i-1$ and $D^r g_{i}(b)=\delta_{r,0}$, \, $r=0,\dots,m-1-i$, with $a$ and $b$ denoting the endpoints of the given interval.
\end{rem}

In view of the fact that Hermite interpolation problems are unisolvent in good for design spline spaces,
the above properties TF-A and TF-B provide a practical way of computing the transition functions.
In particular, from TF-A it is immediate to see that each transition function $f_i$, $i=2,\dots,m+K$, has a nontrivial expression (\ie it is neither the constant function zero or one) in $[t_i,t_{i+m-1}]$. Hence, the restriction of $f_i$ to such interval can be determined by Hermite interpolation, imposing that the left-hand and right-hand pieces of $f_i$ join with continuity $\mu_j$ at each breakpoint $x_j$ contained in {$(t_i,t_{i+m-1})$}, and the endpoint conditions in TF-B. More precisely, let $p_i$ be the index of the break-point associated with a knot $t_i$ such that $t_{i}<t_{i+m-1}$ and let
$x_{p_i},\dots, x_{p_{i+m-1}}$ be the break-points of $\bDelta$ contained in $[t_i,t_{i+m-1}]$.
Then the function $f_{i}$ consists of $p_{i+m-1}-p_i$ pieces and shall satisfy the
continuity conditions and endpoint properties
\begin{equation}\label{eq:cond_tf}
\begin{array}{ll}
D^r f_{i,p_i}(x_{p_i})=0, & r=0,\dots,k_{i}^R,\\ \\
D^r f_{i,j-1}(x_{j})=D^r f_{i,j}(x_{j}), & r=0,\dots,k_{j}, \quad j=p_i+1,\dots,p_{i+m-1}-1, \\ \\
D^r f_{i,p_{i+m-1}-1}(x_{p_{i+m-1}})=\delta_{r,0}, & r=0,\dots,k_{i+m-1}^L,
\end{array}
\end{equation}
where
$f_{i,j}$ is the restriction of $f_{i}$ to the interval $[x_j,x_{j+1}]$, $j=p_i,\dots,p_{i+m-1}-1$,
$k_{i}^R$ and $k_{i+m-1}^L$ are given in \eqref{eq:ki} and $k_{j}\coloneqq m-\mu_{j}-1$. 
{The above conditions \eqref{eq:cond_tf} uniquely determine $f_i$. The proof is analogous to that presented in \cite[Proposition 1]{BCM2017} for the multi-degree polynomial case, therefore it is omitted here.}

\begin{rem}
We would like to draw the readers' attention towards the immediate benefits of computing the B-spline basis from the transition functions as in \eqref{eq:N}. The first of them relates to the fact that the support of the transition functions (intended as the largest interval in which a transition function is neither zero or one) is shorter than the support of the B-spline functions. For instance, in a spline space with $m$-dimensional section spaces and where all knots have multiplicity one, the support of a transition function consists of $m-1$ intervals, whereas the support of a B-spline function of $m$ intervals.
In addition, the fact that $f_1\equiv 1$ guarantees that the basis formed by the functions \eqref{eq:N} enjoys the partition of unity property.
\end{rem}

The computational steps for obtaining $f_i$ can be detailed as follows.
Let $	u_{j,1},u_{j,2},\dots,u_{j,m}$ be a set of generators of the $m$-dimensional EC-space on
 $[x_j,x_{j+1}]$, $j=p_i,\dots,p_{i+m-1}-1$, and $b_{i,s}$, $s=1,\dots,m (p_{i+m-1}-p_i)$ be the coefficients of the local expansion of $f_{i,j}$ with respect to them in such a way that
\begin{equation}\label{eq:expansion}
{f_{i,j}(x)=\sum_{h=1}^m b_{i,m(j-p_i)+h} \, u_{j,h}(x), \quad  x\in [x_j,x_{j+1}], \quad j=p_i,\dots,p_{i+m-1}-1.}
\end{equation}
The coefficients of such expansion can be determined by solving the linear system \eqref{eq:cond_tf}, which, in matrix form reads as
\begin{equation}\label{eq:linear_system}
\matr{A}\vect{b}=\vect{c},
\end{equation}
with
\vspace{-0.1cm}

\begin{widetext}
\begin{equation}\label{eq:Af}
 \matr{A} \coloneqq
\begin{pmatrix}
\matr{A}_{p_i}(x_{p_i}) & & & & \\
\matr{A}_{p_i}(x_{p_i+1}) & -\matr{A}_{p_i+1}(x_{p_i+1}) & & & \\
 & \matr{A}_{p_i+1}(x_{p_i+2}) & -\matr{A}_{p_i+2}(x_{p_i+2}) & & \\
 & & \ddots & \ddots & \\
 & & & \matr{A}_{p_{i+m-1}-2}(x_{p_{i+m-1}-1}) & -\matr{A}_{p_{i+m-1}-1}(x_{p_{i+m-1}-1}) \\
 & & & & \matr{A}_{p_{i+m-1}-1}(x_{p_{i+m-1}}) \\
\end{pmatrix}
\end{equation}
\end{widetext}
$$
\vect{b} \coloneqq \left(b_{i,1},\dots,b_{i,m(p_{i+m-1}-p_i)}\right)^T,
\qquad
\vect{c} \coloneqq (0,\dots,0,1,\underbrace{0,\dots,0}_{k_{i+m-1}^L\text{ times}})^T,
$$
where $\matr{A}_j(x_{h})$, $h=j,j+1$, contains the first $(k_{h}+1)$ rows of the Wronski matrix relative to the $u_{j,1},\dots,u_{j,m}$ evaluated at $x_h$. In particular, it is a matrix of dimensions $(k_{h}+1)\times m$ whose $r$th row, $r=1,\dots,k_{h}+1$, is $\left(D^{r-1} u_{j,1}(x_{h}), \dots, D^{r-1} u_{j,m}(x_{h}) \right)$.
From a practical point of view, it should hence be noted that also the derivatives of $u_{j,1},\dots,u_{j,m}$ must be evaluated for the setup of the linear systems \eqref{eq:linear_system}.

\begin{rem}[Derivatives of a transition function] \label{rem:der}
The solution of the linear system \eqref{eq:linear_system} not only yields the transition function $f_i$, but also its derivatives. In particular, by differentiating the two sides of equation \eqref{eq:expansion}, these can be computed as a linear combination of the derivatives of the generators $u_{j,1},\dots,u_{j,m}$ with coefficients $\vect{b}$ .
\end{rem}

For the sake of clarity, in the following we summarize the described procedures concerning the construction and the evaluation of the transition functions, as well as the evaluation of the Chebyshevian spline, in algorithmic form.


\begin{alg}[Construction of the transition functions]\label{alg:eval_cf}
Let $S(\TTS_m,\MMS,\bDelta)$ be a Chebyshevian spline space and $\bDelta^*$ be the corresponding  extended partition.
To determine the coefficients of the transition functions $f_i$, $i=2,\dots,m+K$, \wrt the EC-space generators $\{u_{j,1},\dots,u_{j,m}\}$, $j=0,\dots,q$ proceed as follows:
\vspace{-0.1cm}
\begin{enumerate}[label=\arabic*.]
\item For $j=0,\dots,q$, calculate the blocks $A_j(x_h)$, $h=j,j+1$ in \eqref{eq:Af}, where $A_j(x_h)$ contains the first $(k_h+1)$ rows of the Wronski matrix relative to 
$u_{j,1}(x),\dots,u_{j,m}(x)$
evaluated at $x_h$.

\vspace{-0.12cm}
\item For each transition function $f_i$, $i=2,\dots,m+K$, assemble the matrix $A$ in  \eqref{eq:Af} and solve the linear system \eqref{eq:linear_system}.
\end{enumerate}
\end{alg}

\vspace{-0.25cm}
{
\begin{alg}[Evaluation of the transition functions]\label{alg:eval_f}
To evaluate at $\bar{x}\in[t_{\ell},t_{\ell+1})=[x_{p_\ell},x_{p_{\ell+1}})$ the (nontrivial) transition functions $f_{\ell-m+2},\dots,f_{\ell}$, it is necessary to evaluate their pieces $f_{\ell-m+2,p_\ell},\dots,f_{\ell,p_\ell}$, proceeding as follows:
\begin{enumerate}[label=\arabic*.]
\vspace{-0.1cm}
\item \textbf{for} $i=\ell-m+2,\dots,\ell$
  \begin{enumerate}
  \vspace{-0.12cm}
  \item[] $\displaystyle f_i(\bar{x}) \leftarrow \sum_{h=1}^m b_{i,m(p_\ell-p_i)+h}\,u_{p_\ell,h}(\bar{x})$ (see eq.\ \eqref{eq:expansion}),
     with $b_{i,m(p_\ell-p_i)+h}$, $h=1,\ldots,m$ computed in Algorithm 1.
  \end{enumerate}
\end{enumerate}
\end{alg}
}


\vspace{-0.2cm}
{
\begin{alg}[Evaluation of a Chebyshevian spline]\label{alg:eval_s}
  Given $\bar{x}\in [a,b]$, to evaluate a spline $s(x)$ of the form  \eqref{eq:spline} (with $N_{i,m}$ satisfying \eqref{eq:N}) at $\bar{x}$:
\vspace{-0.1cm}
\begin{enumerate}[label=\arabic*.]
\item Determine $\ell$ such that $\bar{x}\in [t_\ell,t_{\ell+1})$. \vspace{-0.2cm}
\item Use Algorithm 2 to evaluate the transition functions $f_{\ell-m+2},\dots,f_{\ell}$ at $\bar{x}$. \vspace{-0.2cm}
\item Perform the linear combination, that is $s(\bar x) \leftarrow c_{\ell-m+1}+\sum_{i=\ell-m+2}^\ell (c_i-c_{i-1})f_i(\bar x)$.
\end{enumerate}
\end{alg}}

\vspace{-0.2cm}
{
\begin{alg}[Evaluation of the nontrivial B-splines]\label{alg:eval_s}
Given $\bar{x}\in [a,b]$, to evaluate all the nonzero B-splines $N_{i,m}(\bar{x})$:
\vspace{-0.1cm}
\begin{enumerate}[label=\arabic*.]
\item Determine $\ell$ such that $\bar{x}\in [t_\ell,t_{\ell+1})$. \vspace{-0.2cm}
\item $N_{\ell-m+1,m}(\bar{x}) \leftarrow 1 - \sum_{h=1}^m b_{\ell-m+2,m(p_\ell-p_{\ell-m+2})+h}\,u_{p_\ell,h}(\bar{x})$   \vspace{-0.1cm}
\item \textbf{for} $i=\ell-m+2,\dots,\ell-1$
  \begin{enumerate}
  \vspace{-0.12cm}
    \item[] $\displaystyle N_{i,m}(\bar{x}) \leftarrow \sum_{h=1}^m
      (b_{i,m(p_\ell-p_i)+h} - b_{i+1,m(p_\ell-p_{i+1})+h})\,u_{p_\ell,h}(\bar{x})$ 
  \end{enumerate}
\item $N_{\ell,m}(\bar{x}) \leftarrow \sum_{h=1}^m b_{\ell,h}\,u_{p_\ell,h}(\bar{x})$   \vspace{-0.2cm}
\end{enumerate}
\end{alg}
}

We conclude this section with three instructive examples that aim at presenting in greater detail the computational steps for the construction of the transition functions and of the B-spline basis from them.
The first example is concerned with polynomial splines, being them just special instances of Chebyshevian splines.

\begin{exmp}[Polynomial B-splines]
Let us consider the space of polynomial splines $S(\PP_m,\MMS,\bDelta)$, where each piece is spanned by
the space $\PP_m$ of polynomials of degree $m-1$ (\ie, of degree at most $m-1$) and where, for the sake of simplicity, every break-point in $\bDelta$ has multiplicity $\mu_i=1$. The transition functions $f_{1}\equiv 1$ and $f_{i}$, $i=2,\dots,m+K$, are the degree-$(m-1)$ piecewise polynomial functions obtained as the unique solution of the linear systems
\begin{equation}\label{eq:cont_cond}
\begin{array}{ll}
D^r f_{i,i}(t_i)=0, & r=0,\dots,m-2,\\
D^r f_{i,j-1}(t_{j})=D^r f_{i,j}(t_{j}), & r=0,\dots,m-2, \quad j=i+1,\dots,i+m-2,\\
D^r f_{i,i+m-2}(t_{i+m-1})=\delta_{r,0}, & r=0,\dots,m-2.
\end{array}
\end{equation}

Once solved the above systems, the B-spline basis functions $N_{i,m}$ are given by \eqref{eq:N}.

For example, take $[a,b]=[0,3]$, $m=4$ and a clamped uniform extended partition $\bDelta^*=\{0,0,0,0,1,2,3,3,3,3\}$. The considered space has dimension 6 and associated transition functions $f_1\equiv 1$ and $f_i$, $i=2,\dots,6$ each of which is the solution of a linear system of the form \eqref{eq:linear_system}, where the entries of each block of the coefficient matrix $A$ depend on the generators of the corresponding section space on the nontrivial knot intervals $[t_j,t_{j+1}]\subseteq[t_i,t_{i+m-1}]$. Each block $A_j(x_h)$ has as many rows as the number of continuity conditions that $f_i$ must satisfy at $x_h$, and as many columns as the dimension of the section spaces, that is $4$.  If the  generators of the section space on $[t_j,t_{j+1}]$ are taken to be $1$, $x-t_j$, $(x-t_j)^2$, $(x-t_j)^3$, then the matrices obtained for the computation of $f_2$ and $f_3$ (respectively the matrices $A$ to the left-hand side and right-hand side below)
read as
\begin{equation}
\newcommand*{\tempz}{\multicolumn{1}{|c}{0}}
\newcommand*{\tempu}{\multicolumn{1}{|c}{1}}
\newcommand*{\tempmu}{\multicolumn{1}{|c}{-1}}
\newcommand*{\tempzl}{\multicolumn{1}{c|}{0}}
\newcommand*{\tempul}{\multicolumn{1}{c|}{1}}
\newcommand*{\temptl}{\multicolumn{1}{c|}{3}}
\newcommand*{\tempsl}{\multicolumn{1}{c|}{6}}
A=\left[\,
\begin{array}{cccc}
\cline{1-4}
     \tempu  &   0  &   0  &   \tempzl\\\cline{1-4}
     \tempu   &  1  &   1  &   \tempul\\
     \tempz   &  1  &   2  &   \temptl\\
     \tempz   &  0  &   2  &   \tempsl\\\cline{1-4}
\end{array}
\,\right]
\qquad
\hbox{and}
\qquad
A=\left[\,
\begin{array}{cccccccc}
\cline{1-4}
     \tempu  &   0  &   0  &   0  &   \tempz  &   0  &   0  &   0\\
     \tempz  &   1   &  0   &  0   &  \tempz   &  0   &  0  &  0\\\cline{1-8}
     \tempu   &  1   &  1   &  1  &  \tempmu   &  0   &  0   &  \tempzl\\
     \tempz   &  1   &  2   &  3   &  \tempz  &  -1   &  0   &  \tempzl\\
     \tempz   &  0   &  2   &  6   &  \tempz  &   0   & -2   &  \tempzl\\\cline{1-8}
     0   &  0   &  0   &  0  &   \tempu   &  1   &  1   &  \tempul\\
     0   &  0   &  0   &  0  &   \tempz   &  1   &  2   &  \temptl\\
     0   &  0   &  0   &  0   &  \tempz   &  0   &  2   &  \tempsl  \\ \cline{5-8}
\end{array}\,
\right],
\end{equation}
whereas for the computation of $f_4$ we obtain
\begin{equation}
\newcommand*{\tempz}{\multicolumn{1}{|c}{0}}
\newcommand*{\tempu}{\multicolumn{1}{|c}{1}}
\newcommand*{\tempmu}{\multicolumn{1}{|c}{-1}}
\newcommand*{\tempzl}{\multicolumn{1}{c|}{0}}
\newcommand*{\tempul}{\multicolumn{1}{c|}{1}}
\newcommand*{\temptl}{\multicolumn{1}{c|}{3}}
\newcommand*{\tempsl}{\multicolumn{1}{c|}{6}}
A=\left[\,
\begin{array}{cccccccccccc}
\cline{1-4}
     \tempu   &  0  &   0   &  0   &  \tempz    &  0   &  0   &  0   &  0   &  0   &  0  &   0\\
     \tempz   &  1  &   0  &   0   &  \tempz    &  0  &   0  &   0   &  0  &   0   &  0   &  0\\
     \tempz   &  0  &   2   &  0   &  \tempz   &   0  &   0  &   0  &   0  &   0  &   0  &   0\\\cline{1-8}
     \tempu   &  1  &   1  &   1 &   \tempmu  &   0  &   0   &  0  &   \tempz  &   0   &  0  &   0\\
     \tempz   &  1  &   2  &   3   &  \tempz  &  -1 &    0  &   0   &  \tempz   &  0   &  0   &  0\\
     \tempz   &  0 &    2  &   6  &   \tempz  &   0  &  -2  &   0  &   \tempz  &   0   &  0  &   0\\\cline{1-12}
     0   &  0  &   0 &    0   &  \tempu  &   1  &   1  &   1  &  \tempmu  &   0  &   0  &   \tempzl\\
     0   &  0   &  0  &   0   &  \tempz  &   1   &  2   &  3  &   \tempz  &  -1  &   0    & \tempzl\\
     0   &  0   &  0  &   0   &  \tempz  &   0   &  2   &  6  &   \tempz   &  0   & -2    & \tempzl\\\cline{5-12}
     0   &  0   &  0   &  0   &  0  &   0   &  0   &  0  &   \tempu   &  1   &  1   &  \tempul\\
     0    & 0   &  0   &  0   &  0   &  0   &  0   &  0   &  \tempz  &   1   &  2   &  \temptl\\
     0    & 0   &  0   &  0   &  0   &  0   &  0   &  0   &  \tempz  &  0   &  2  &   \tempsl \\\cline{9-12}
\end{array}
\right].
\end{equation}
For the sake of conciseness we omit the expression of the matrices for
$f_5$ and $f_6$, which have the same size and analogous structure as those for $f_3$ and $f_2$, respectively.
\end{exmp}

The next two examples deal with nontrivial instances of Chebyshevian splines. Precisely, in Example \ref{ex:Bspline} we
illustrate the construction of a single B-spline function for a Chebyshevian spline having polynomial, trigonometric and hyperbolic pieces, whereas in Example \ref{ex:mu_zero} we are concerned with the construction of the B-spline basis for a Chebyshevian spline having section spaces spanned by rational functions with tension parameters, and breakpoints with multiplicity equal to zero.

\begin{figure*}[h!]
\centering
\subfigure
{\includegraphics[width=0.8\textwidth/3]{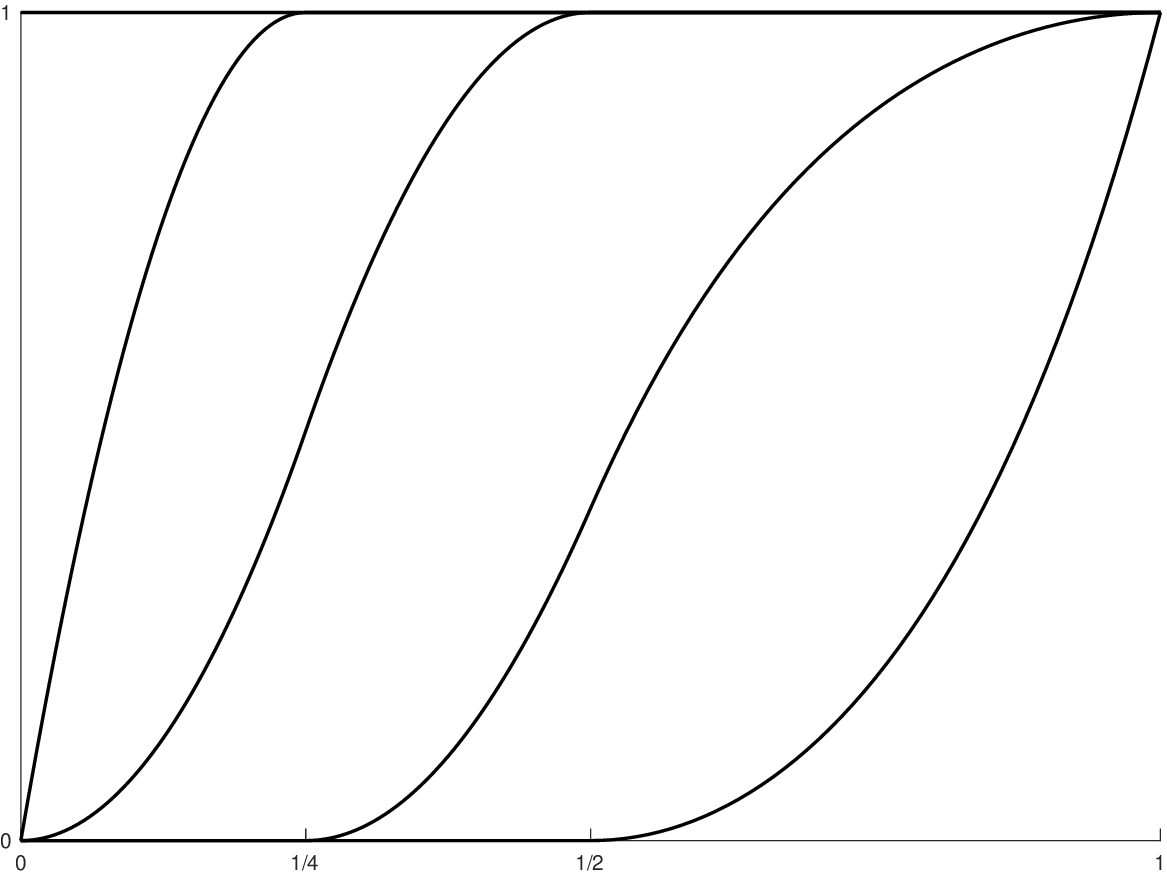}\label{fig:fN_f}}
\hfill
\subfigure
{\includegraphics[width=0.8\textwidth/3]{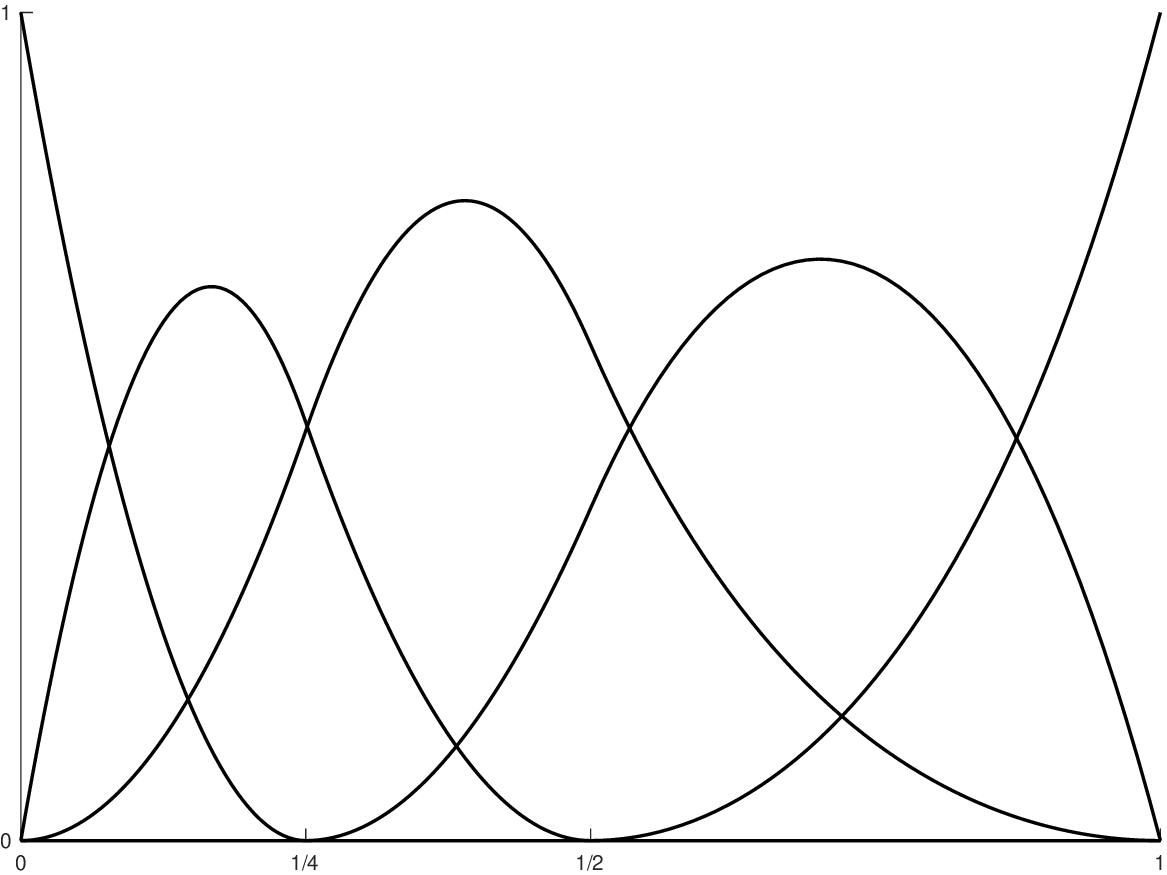}\label{fig:fN_N}}
\hfill
\subfigure
{\includegraphics[width=0.8\textwidth/3]{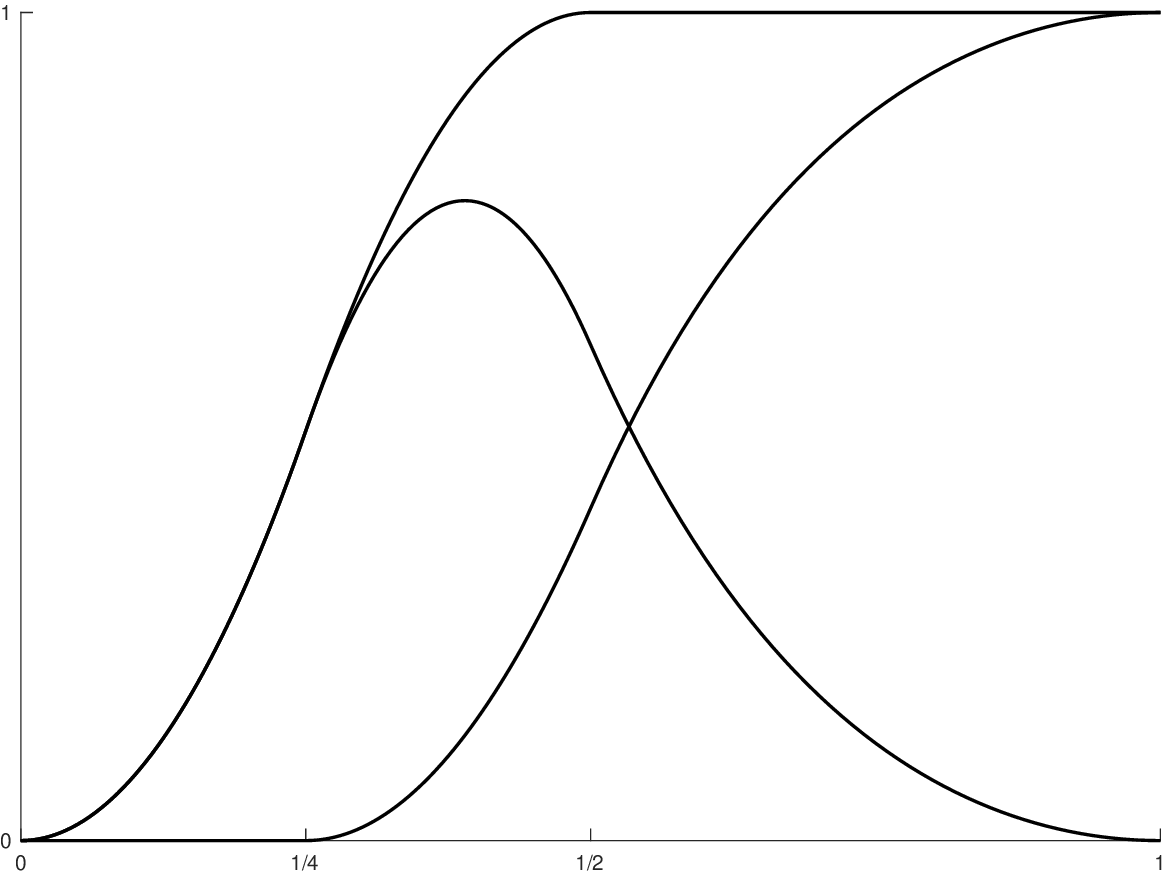}\label{fig:fN_fNpart}}
\caption{
Transition functions (left) and B-spline basis (center) for the spline space $S(\TTS_3,\MMS,\bDelta)$ of Example \ref{ex:Bspline} with
$\theta=2$ and $\phi=4$. The single B-spline function $N_{3,3}$ and the associated transition functions $f_{3}, f_{4}$ (right).}
\label{fig:fN}
\end{figure*}

\begin{exmp}\label{ex:Bspline}
We consider the ordered sequence $\TTS_3=\{\TT_{0,3},\TT_{1,3},\TT_{2,3}\}$ made of the 3-dimensional EC-spaces
\begin{equation}
\begin{array}{ll}
\TT_{0,3} = \Span\{1,t,t^2\}, & t=t(x) \coloneqq x-x_{0}, \quad x \in [x_0,x_{1}], \\ \\
\TT_{1,3} = \Span\{1,\cos(\theta t),\sin(\theta t)\}, & t=t(x) \coloneqq
  x-x_{1}, \quad x \in [x_1,x_{2}], \;\theta(x_2-x_1)<L=\pi, \\ \\
\TT_{2,3} = \Span\{1,\cosh(\phi t),\sinh(\phi t)\}, & t=t(x) \coloneqq x-x_{2},
  \quad x \in [x_2,x_{3}],
\end{array}
\end{equation}
and, for the ease of presentation, we assume all break-points to have multiplicity 1. Being $\bDelta^*=\{t_1,\ldots,t_8\}$ the extended partition, we derive the expression of the B-spline function $N_{3,3}=f_{3}-f_{4}$ having support $[t_{3},t_{6}]=[x_{0},x_{3}]$.
Expanding $f_{j}$, $j=3,4$, w.r.t. the EC-space generators $\{u_{\ell,1},u_{\ell,2},u_{\ell,3}\}$, $\ell=0,1,2$, we obtain
\[
f_{j}(x) =
\begin{cases}
0, & x \leq t_{j}, \\
\displaystyle \sum_{k=1}^{3} b_{j,k} \, u_{j-3,k}(t(x)), & t_{j} \leq x <t_{j+1}, \\
\displaystyle \sum_{k=1}^{3} b_{j,k+3} \, u_{j-2,k}(t(x)), & t_{j+1} \leq x <t_{j+2}, \\
1, & x \geq t_{j+2}.
\end{cases}
\]
Next we compute the coefficients $b_{j,k}$ of the transition function $f_{j}$, $j=3,4$, by solving the two $6\times6$ linear systems of the form \eqref{eq:linear_system} that arise from the endpoint and continuity conditions
\begin{equation}\label{system_66}
\begin{array}{llll}
& f_{j,1}(x_{j-3})=0, &\qquad & Df_{j,1}(x_{j-3})=0,\\[1ex]
& f_{j,1}(x_{j-2})=f_{j,2}(x_{j-2}), &\qquad& Df_{j,1}(x_{j-2})=Df_{j,2}(x_{j-2}), \\[1ex]
& f_{j,2}(x_{j-1})=1, &\qquad& Df_{j,2}(x_{j-1})=0.
\end{array}
\end{equation}

\noindent
Solving \eqref{system_66} for $f_{3}$ and  $f_{4}$, we respectively obtain
\begin{equation}
\begin{array}{l}
b_{3,1} = b_{3,2} = 0, \\
b_{3,3} = {\displaystyle\frac{1}{h_{0}^2+2\frac{h_{0}}{\theta}\tan\left(\frac{\theta}{2}h_{1}\right)}}, \\
b_{3,4} = 1+{\displaystyle\frac{1}{\cos(\theta h_{1})-\frac{\theta}{2}h_{0}\sin(\theta h_{1})-1}}, \\
b_{3,5} = (1-b_{3,4})\cos(\theta h_{1}), \\
b_{3,6} = {\displaystyle\frac{1}{\frac{\theta}{2}h_{0}+\tan\left(\frac{\theta}{2}h_{1}\right)}}, \\
\end{array}
\qquad
\hbox{and}
\qquad
\begin{array}{l}
b_{4,1} = -b_{4,2} = {\displaystyle\frac{1}{1- \cos(\theta h_{1}) +\frac{\theta}{\phi} \sin(\theta h_{1}) \tanh\left(\frac{\phi}{2}h_{2}\right)}}, \\
b_{4,3} = 0, \\
b_{4,4} = 1+{\displaystyle\frac{1}{\cosh(\phi h_{2}) +\frac{\phi}{\theta} \tan\left(\frac{\theta}{2}h_{1}\right) \sinh(\phi h_{2}) -1}}, \\
b_{4,5} = (1-b_{4,4})\cosh(\phi h_{2}), \\
b_{4,6} = {\displaystyle\frac{1}{\frac{\phi}{\theta} \tan\left(\frac{\theta}{2}h_{1}\right) +\tanh\left(\frac{\phi}{2}h_{2}\right)}}, \\
\end{array}
\end{equation}
where $h_i=t_{i+4}-t_{i+3}$, $i=0,1,2$. An illustration of the transition functions and the related B-spline basis for the space $S(\TTS_3,\MMS,\bDelta)$ with extended partition $\bDelta^*=\{0,0,0,\frac{1}{4},\frac{1}{2},1,1,1\}$ and parameters $\theta=2, \phi=4$, is given in Figure \ref{fig:fN}.
\end{exmp}

\begin{figure*}[h!]
\centering
\subfigure[]{
{\includegraphics[width=0.8\textwidth/3]{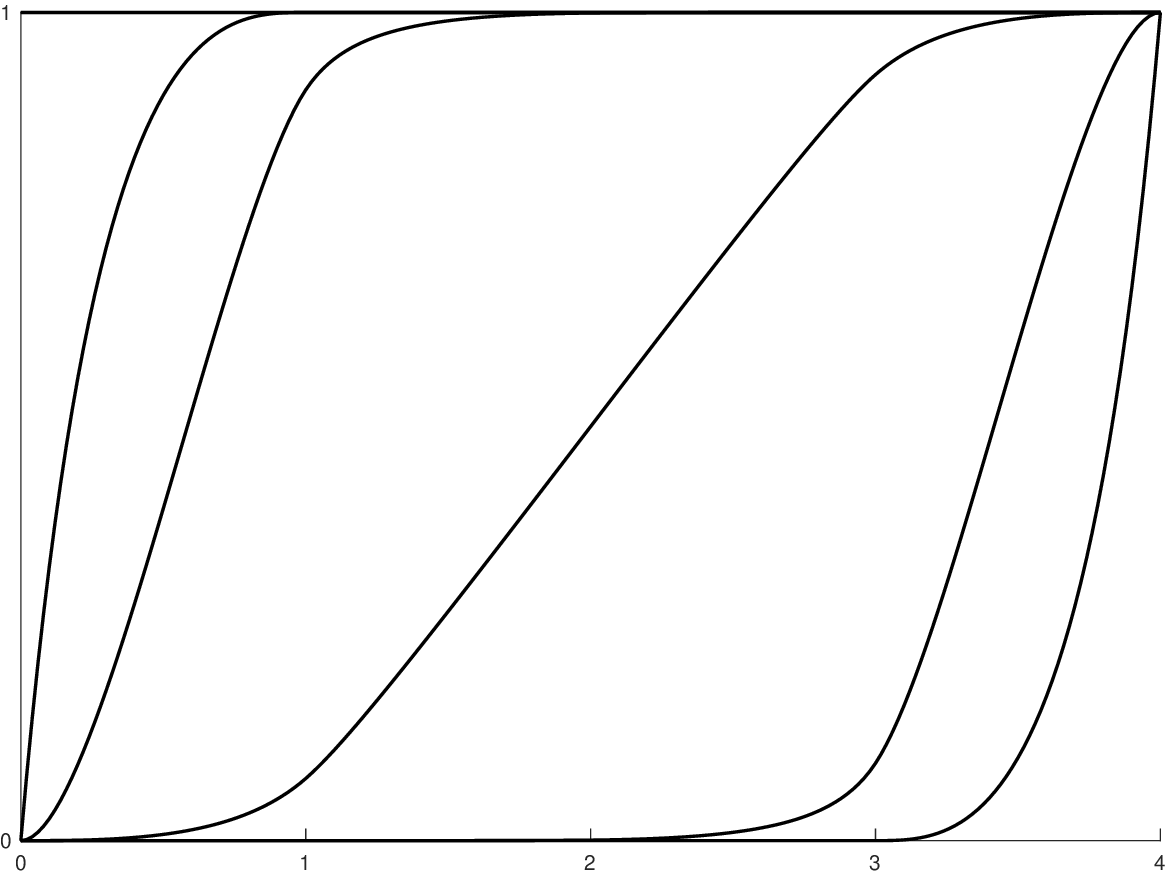}\label{fig:B-spline_zero_knot-crop}}
\hfill
}
\hfill
\subfigure[]{
{\includegraphics[width=0.8\textwidth/3]{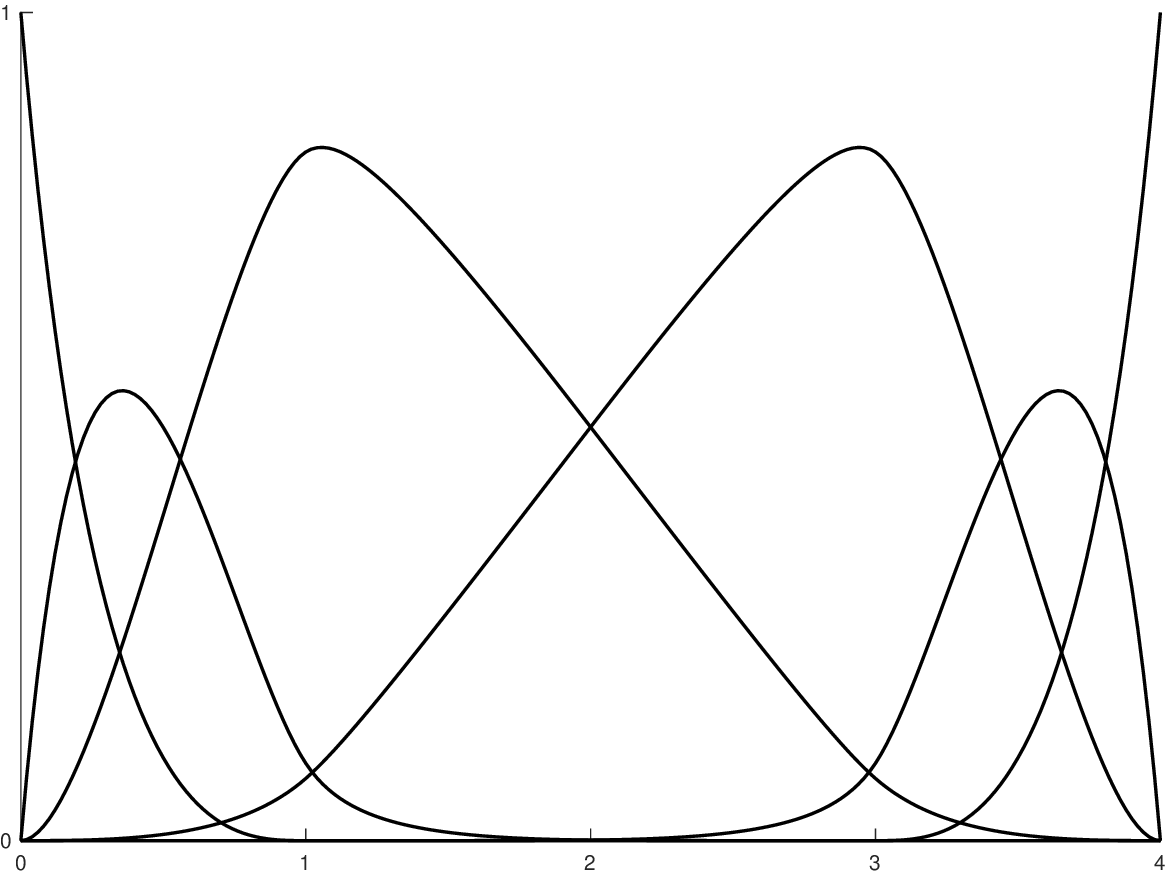}\label{fig:B-spline_zero_knot-crop}}
\hfill
}
\hfill
\subfigure[]{
{\includegraphics[width=0.8\textwidth/3]{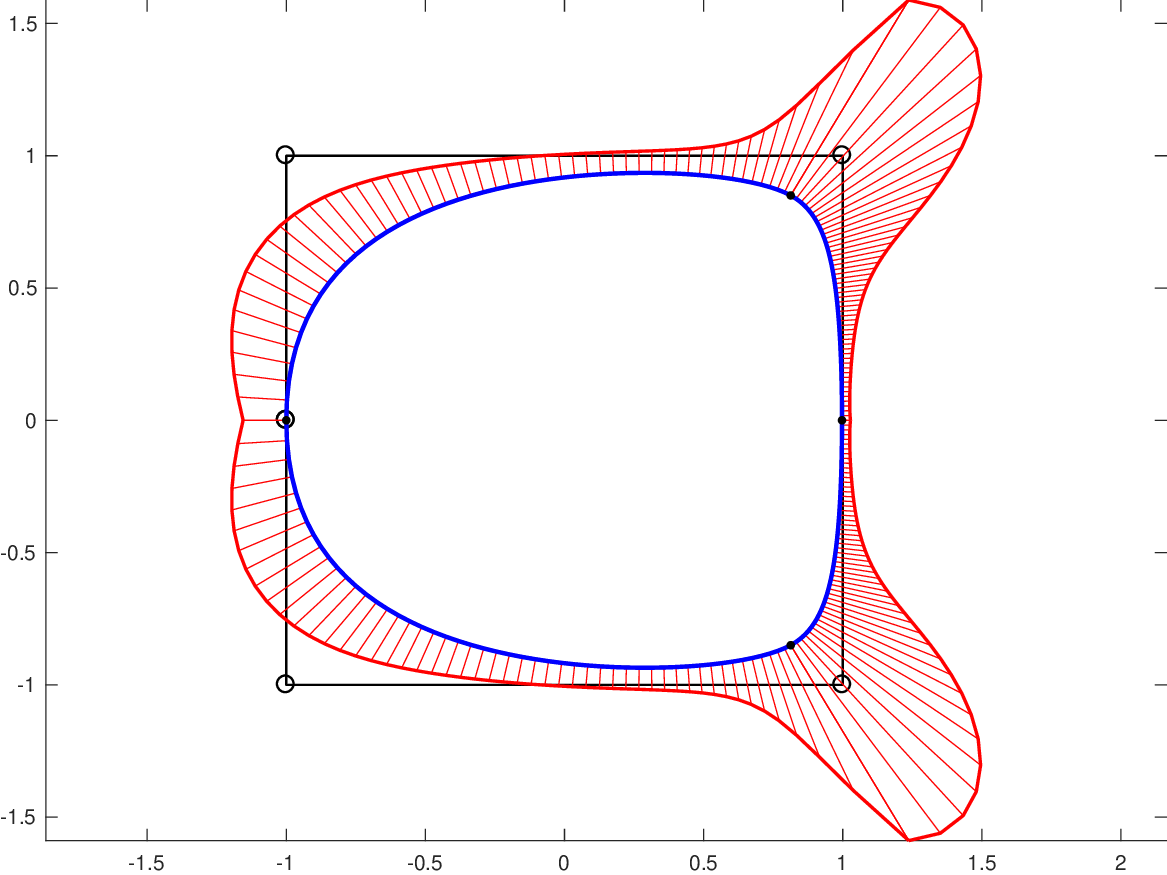}\label{fig:zero_knot-crop}}
}
\caption{
Transition functions (left) and corresponding B-spline functions (center) for the spline space $S(\TTS_4,\MMS,\bDelta)$ of Example \ref{ex:mu_zero} with $\nu_0=\nu_3=4$ and $\nu_1=\nu_2=6$. The closed spline curve obtained from the highlighted square control polygon, and its curvature comb (right).
}
\label{fig:zero_knot}
\end{figure*}

\begin{exmp}\label{ex:mu_zero}
When defining Chebyshevian splines it is possible (and in fact usual) to allow the multiplicity $\mu_i$ of a knot $x_i$ to be equal to zero. In fact, when there are two different $m$-dimensional EC-spaces on $[x_{i-1},x_{i}]$ and $[x_i,x_{i+1}]$, we can have a piecewise function with continuity $C^{m-1}$ at $x_i$.
The break-point $x_i$ is said to have \emph{zero multiplicity} because there is no corresponding knot $t_j$ in the extended partition $\bDelta^*$.\\
To illustrate such situation, consider an ordered sequence $\TTS_4$ of 4-dimensional EC-spaces of the form
$$
\TT_{i,4} = \Span \left\{1, t, \frac{(1-t)^{3}}{1+(\nu_{i}-3)(1-t)t}, \frac{t^{3}}{1+(\nu_{i}-3)(1-t)t}\right\},
\qquad t=\frac{x-x_i}{x_{i+1}-x_i}, \quad x \in [x_i, x_{i+1}], \quad \nu_i \geq 3.
$$
The parameters $\nu_i$ are known to have a \qtext{tension} effect and therefore can be used to efficiently solve shape-preserving interpolation problems \cite{Costantini2010592}.
By setting $[a,b]=[0,4]$, $\bDelta=\{1,2,3\}$ and $\MMS=(1,0,1)$, the ordered sequence $\TTS_4$ contains four different spaces $\TT_{i,4}$, $i=0,\ldots,3$ of the previously considered form.
Figure \ref{fig:zero_knot}, obtained by selecting the tension parameters $\nu_0=\nu_3=4$ and $\nu_1=\nu_2=6$, shows the transition functions of the spline space $S(\TTS_4,\MMS,\bDelta)$, the associated B-spline basis and
a closed spline curve generated with $6$ control points on a square.
The curvature comb emphasizes that the curve is $C^2$ at the break-points
$1$ and $3$, and $C^3$ at the break-point $2$, which has multiplicity zero.
\end{exmp}

\section{Knot insertion and order elevation via transition functions} \label{sec:knotins_ordelev}
In this section we illustrate how to use the transition functions to efficiently perform knot insertion and dimension elevation in Chebyshevian spline spaces good for design.

\subsection{Knot insertion} \label{sec:knot_ins}
Starting from a spline space good for design, we can insert one knot either inside a breakpoint interval or in correspondence of another existing knot (thus increasing the multiplicity within the limit $\mu_i< m$) generating another spline space good for design.
In particular, denoted by $S(\TTS_m,\MMS,\bDelta)$ the initial space, with
$\bDelta^*=\{t_i\}_{i=1,\dots,2m+K}$, by $\hat{t}$ the inserted knot such that $\hat{t}\in [t_{\ell},t_{\ell+1})$, by $S(\rwh{\TTS}_m,\rwh{\MMS},\rwh{\bDelta})$ the space obtained after knot insertion, with
$\rwh{\bDelta}^*=\{\hat{t}_i\}_{i=1,\dots,2m+K+1}$, and by $r$ , $1 \leq r < m$ the multiplicity of $\hat t \in \hat \Delta^*$, the respective B-spline bases
$N_{j,m}$ and $\rwh{N}_{j,m}$ satisfy the following relation
\begin{equation}\label{eq:ki_basis}
\begin{array}{l}
N_{i,m}(x) = \alpha_{i,m}\,\rwh{N}_{i,m}(x)
  + (1-\alpha_{i+1,m})\,\rwh{N}_{i+1,m}(x), \quad i=1,\dots,m+K,\\
\end{array}
\end{equation}
with coefficients $0\leq \alpha_{i,m} \leq 1$, $i=1,\dots,m+K+1$.
For a spline function
$$
s(x) = \sum_{i=1}^{m+K} c_i\,N_{i,m}(x) = \sum_{i=1}^{m+K+1} \hat{c}_i\,\rwh{N}_{i,m}(x),
$$
an immediate consequence of \eqref{eq:ki_basis} is that the coefficients of the two expansions are given by
\begin{equation}\label{eq:ki_spl_coeff}
\hat{c}_i =
\begin{cases}
c_i, & i \leq \ell-m+1,\\
\alpha_{i,m}\,c_i + (1-\alpha_{i,m})\,c_{i-1}, & \ell-m+2 \leq i \leq \ell-r+1,\\
c_{i-1}, & i \geq \ell-r+2.
\end{cases}
\end{equation}

The expressions of the coefficients $\alpha_{1,m},\dots,\alpha_{m+K+1,m}$, whose computation is in general nontrivial, can however efficiently be determined in dependence of the transition functions. 
From \eqref{eq:N} and \eqref{eq:ki_basis} one gets the expression
\begin{equation} \label{eq:f_i_knot_ins}
f_{i}-f_{i+1}=\alpha_{i,m}(\hat f_{i}- \hat f_{i+1})+(1-\alpha_{i+1,m})(\hat f_{i+1}- \hat f_{i+2}),
\end{equation}
which after differentiating $k_{i}^R+1$ times, with $k_{i}^R$ defined in \eqref{eq:ki}, and evaluating the result at $t_i$, yields
\begin{equation}\label{eq:ki_coef}
\alpha_{i,m}=\left\{
\begin{array}{ll}
1, & i \leq \ell-m+1, \\
{\displaystyle\frac{D_+^{k_{i}^R+1}f_{i,p_i}({t_{i}})}{D_+^{k_{i}^R+1}\hat{f}_{i,p_i}({t_{i}})}}, & \ell-m+2 \leq i \leq \ell-r+1, \\
0, & i \geq \ell-r+2,
\end{array} \right.
\end{equation}
where
$f_{i,p_i}$ and $\hat{f}_{i,p_i}$ are the first non-trivial pieces of the transition functions in  \eqref{eq:N} which determine $N_{i,m}$ and $\rwh{N}_{i,m}$.

\begin{rem}
The knot-insertion coefficients $\alpha_{i,m}$ in \eqref{eq:ki_coef} are determined as the ratio of derivatives of the transition functions, which requires that these derivatives be nonzero at the evaluation point $t_i$. Observing that, at $t_i$, both $f_i$ and $\hat{f}_i$ have the same order of continuity $k_{i}^R$, there follows that we shall take their derivatives of order $k_{i}^R+1$.
The right-hand side evaluation is thus explained by the fact that the considered functions are continuous up to order $k_{i}^R$ only at the evaluation point $t_i$.
\end{rem}

\begin{rem}
An equivalent way for computing the coefficients $\alpha_{i,m}$, $\ell-m+2 \leq i \leq \ell-r+1$ in \eqref{eq:ki_coef} is:
\begin{equation}
\alpha_{i,m}={\displaystyle\frac{D_-^{k_{i+m-1}^L+1}f_{i,p_{i+m-1}-1}({t_{i+m-1}})}{D_-^{k_{i+m-1}^L+1}\hat{f}_{i+1,p_{i+m-1}}({t_{i+m-1}})}}.
\label{eq:ki_coef2}
\end{equation}

Formula \eqref{eq:ki_coef2} is useful when the additional knots in the extended partition are noncoincident.
In this case, we should use the above formula when inserting a knot $\hat{t} \equiv b$ (in order to switch from an extended partition with distinct additional knots to one with coincident additional knots) or $\hat{t} > b$.
\end{rem}

The steps to perform knot insertion can be summarized as follows.

\begin{alg}[Knot insertion]\label{alg:knot_ins}
Let $S(\TTS_m,\MMS,\bDelta)$ be a given Chebyshevian spline space, $\bDelta^*$ an extended partition, and $s$ a spline function of the form \eqref{eq:spline}.
Given $\hat{t}\in[a,b]$, to determine the coefficients $\hat{c}_i$, $i=1,\dots,m+K+1$, obtained by inserting $\hat{t}$ in $\bDelta^*$ once, proceed as follows:
\begin{enumerate}[label=\arabic*.]
\item determine $\ell$ such that $\hat{t}\in [t_{\ell},t_{\ell+1})$;
\item compute the coefficients $\alpha_i$, $i=\ell-m+2,\dots,\ell-r+1$ through \eqref{eq:ki_coef};
\item use \eqref{eq:ki_spl_coeff} to compute the coefficients $\hat{c}_i$, $i=1,\dots,m+K+1$.
\end{enumerate}
\end{alg}

To conclude this section we would like to draw the readers' attention to the computational efficiency of the transition functions.
In fact the local nature of both the B-spline basis and the transition functions makes so that
only a few of these functions change after knot insertion and the local nature of the proposed method (in which each transition function is computed individually and independently of the others) makes it necessary to recalculate only the ones which need to be updated. Moreover the evaluation of the derivatives of transition functions appearing in \eqref{eq:ki_coef}
is immediate as already observed in Remark \ref{rem:der}.

\subsection{Dimension elevation} \label{sec:ord_elev}
The aim of dimension elevation is to write a spline in a space $S(\TTS_m,\MMS,\bDelta)$, with $m$-dimensional sections, as an element of a space $S(\TTS_{m+r},\widetilde{\MMS},\bDelta)$, with $(m+r)$-dimensional sections ($r\geq1$) and such that
 $S(\TTS_m,\MMS,\bDelta)$ $\subset S(\TTS_{m+r},\widetilde{\MMS},\bDelta)$.
 To do so, it is necessary to assume that both spaces are good for design and thus have a B-spline basis. Unlike for polynomials, however, in general this does not guarantee that dimension elevation will be \qtext{corner-cutting}  (see counter examples in \cite{BCM2020b}) or, in other words, that the control polygon of the higher dimensional spline representation can be obtained cutting the corners of the control polygon of the lower dimensional one.

The following three steps allow for expressing the B-splines $\{N_{i,m}\}_i$ in $S(\TTS_m,\MMS,\bDelta)$ as a linear combination of the B-splines $\{N_{i,m+r}\}_i$ in $S(\TTS_{m+r},\widetilde{\MMS},\bDelta)$, where $\widetilde{\MMS}$ is obtained from $\MMS$ by simply increasing by $r$ the multiplicity of each break-point:
\begin{enumerate}
\item[DE-1] \label{enum:oe_step1} { apply repeated knot insertions to} subdivide the spline into its pieces, each of which is represented in the Bernstein basis;
\item[DE-2] \label{enum:oe_step2} elevate the dimension of each spline piece;
\item[DE-3] \label{enum:oe_step3} remove { the inserted knots up to restoring the original continuity}.
\end{enumerate}

For step DE-1, it is sufficient to recall that a representation in the Bernstein basis can be derived by performing
repeated knot insertions, as illustrated in Section \ref{sec:knot_ins}. { After elevating the dimension of each piece, as per step DE-2, a suitable number of the inserted knots will need to be removed so as to restore the original continuity, as prescribed in DE-3.}
Step DE-2 { is the core of the dimension elevation process, as it }requires relating the Bernstein bases of the section spaces of dimension $m$ and
$m+r$. 

{ To describe step DE-2, f}or the moment we focus on $r=1$ and denote by  $\TT_{n+1}$ and $\TT_{n+2}$ two EC-spaces on $[a,b]$ of dimension $n+1$ and $n+2$, respectively, containing constants and such that $\TT_{n+1}\subset\TT_{n+2}$. In this case dimension elevation is a  corner cutting algorithm and
the respective Bernstein bases $\{B_{i,n}\}_{i=0,\dots,n}$ and $\{B_{i,n+1}\}_{i=0,\dots,n+1}$ satisfy the relation
\begin{equation}\label{eq:order_elev_1_basis}
\begin{array}{l}
B_{i,n}(x) = \gamma_i\,B_{i,n+1}(x) + (1-\gamma_{i+1})\,B_{i+1,n+1}(x), \; \;i=0,\dots,n,\;\; x\in [a,b], \\
\end{array}
\end{equation}
with coefficients $\gamma_0=1$, $\gamma_{n+1}=0$ and $\gamma_i\in {(0,1)} $, $i=1,\dots,n$, see \cite{BCM2020b} and references therein. As a consequence of the above identity, the expansions of a function $p$ in the two spaces
\[
p(x) = \sum_{i=0}^{n} c_i\,B_{i,n}(x) = \sum_{i=0}^{n+1} \tilde{c}_i\,B_{i,n+1}(x), \quad x\in[a,b],
\]
have coefficients
\begin{equation} \label{eq:coef_elev_1}
\tilde{c}_0=c_0,\qquad \tilde{c}_i =(1-\gamma_{i})\,c_{i-1} + \gamma_{i}\,c_i, \quad i=1,\dots,n, \qquad \tilde{c}_{n+1}=c_{n}.
\end{equation}

The dimension elevation coefficients $\gamma_i${, $i = 1,\dots,n$} can conveniently be computed from the transition functions $g_i$ 
 in $\TT_{n+1}$ and
$\tilde{g}_i$  
 in $\TT_{n+2}${, $i=1,\dots,n$}
(see Remark \ref{rem:Bernst} for the notation).
In particular, { differentiating \eqref{eq:order_elev_1_basis} $i$ times and evaluating at $x=a$ yields
$
\gamma_i = \displaystyle \frac{D^{i}B_{i,n}(a)}{D^{i}B_{i,n+1}(a)}, \, i=1,\dots,n, 
$
}
from which we obtain the following expressions of $\gamma_i$ in dependence of the transition functions:
\begin{equation}\label{eq:order_elev_1}
\gamma_i=\frac{D^{i}g_i(a)}{D^{i}\tilde{g}_i(a)}, \quad i=1,\dots,n.
\end{equation}

To raise the section spaces of more than one dimension, that is considering $r>1$, one may perform several successive elevations of one dimension as described above.
However, when dealing with certain Chebyshevian spline spaces, it is convenient to elevate by two dimensions at the same time
(think, e.g., of a situation in which we may want to add the pair of functions $\sin,\cos$ or $\sinh,\cosh$).
To this aim,
let $\TT_{n+1}$ and $\TT_{n+3}$ be respectively $(n+1)$- and $(n+3)$-dimensional EC-spaces on $[a,b]$, such that $\TT_{n+1} \subset \TT_{n+3}$ and let
$\{B_{i,n}\}_{i=0,\dots,n}$ and $\{B_{i,n+2}\}_{i=0,\dots,n+2}$ be their Bernstein bases.
By recalling that $B_{i,n}$ vanishes $i$ times at $a$ and $n-i$ times at $b$,
one can see that the relation between the bases of dimensions $n+1$ and $n+3$ has the form
\begin{equation}\label{eq:de222}
\begin{array}{l}
B_{i,n}(x) = \gamma_i\,B_{i,n+2}(x) + \delta_{i+1}\,B_{i+1,n+2}(x) + \varepsilon_{i+2}\,B_{i+2,n+2}(x), \qquad x\in[a,b].
\end{array}
\end{equation}
{
Since
\begin{equation}
1 \equiv \sum_{i=0}^n B_{i,n}(x)
  = \sum_{i=0}^n \left( \gamma_i\,B_{i,n+2}(x) + \delta_{i+1}\,B_{i+1,n+2}(x) + \varepsilon_{i+2}\,B_{i+2,n+2}(x) \right)
  = \sum_{i=0}^{n+2} B_{i,n+2}(x),
\end{equation}
for all $x\in[a,b]$, we have $\gamma_0=\varepsilon_{n+2}=1$, $\gamma_1+\delta_1=1$, $\delta_{n+1}+\varepsilon_{n+1}=1$, and $\gamma_i+\delta_i+\varepsilon_i=1$, for $i=0,\dots,n+2$.
Moreover, like before, the dimension elevation coefficients $\gamma_i$, $i=1,\dots,n$ and $\delta_i$, $i=2,\dots,n+1$ can conveniently be computed from the transition functions $g_i$, $i=1,\dots,n$, in $\TT_{n+1}$ and
$\tilde{g}_i$, $i=1,\dots,n+1$ in $\TT_{n+3}$.
In particular 
\begin{equation}
\gamma_i =
\displaystyle
\frac{D^{i}B_{i,n}(a)}{D^{i}B_{i,n+2}(a)}, \; i=1,\dots,n,
\qquad \quad 
\delta_{i+1} =
\displaystyle
\frac{D^{i+1}B_{i,n}(a)-\gamma_i D^{i+1}B_{i,n+2}(a)}{D^{i+1}B_{i+1,n+2}(a)}, \; i=1,\dots,n,
\end{equation}
from which we can get the following expressions of the 2-dimension elevation coefficients
\begin{equation}\label{eq:order_elev_2a}
\gamma_i =
\displaystyle
\frac{D^{i}g_i(a)}{D^{i}\tilde{g}_{i}(a)}, \; i=1,\dots,n, 
\qquad \qquad 
\delta_{i+1} =
\displaystyle
\gamma_i-\gamma_{i+1}+\frac{D^{i+1}g_i(a)-\gamma_iD^{i+1}\tilde{g}_{i}(a)}{D^{i+1}\tilde{g}_{i+1}(a)},
\; i=1,\dots,n,
\end{equation}
in dependence of the transition functions.
} 



As an immediate consequence of \eqref{eq:de222}, we get that the expansions
\[
p(x) = \sum_{i=0}^{n} c_i\,B_{i,n}(x) = \sum_{i=0}^{n+2} \tilde{c}_i\,B_{i,n+2}(x), \quad x\in[a,b],
\]
have coefficients such that
\begin{equation}\label{eq:coef_elev_2}
\tilde{c}_i =
\begin{cases}
c_0, & i=0,\\
\gamma_{1}\,c_1 + (1-\gamma_{1})\,c_{0}, & i=1,\\
\gamma_{i}\,c_i + \delta_{i}\,c_{i-1}+(1-\gamma_{i}-\delta_{i})\,c_{i-2}, & i=2,\dots,n,\\
\delta_{n+1}\,c_{n}+(1-\delta_{n+1})\,c_{n-1}, & i=n+1,\\
c_n, & i=n+2.
\end{cases}
\end{equation}
Note that, while from the endpoint properties of the Bernstein bases it is easy to see that the coefficients $\gamma_i$ are nonnegative, there is no condition guaranteeing that also $\delta_i$ and $\varepsilon_i$ be nonnegative too. In particular, as previously recalled, the assumption that both spaces $\TT_{n+1}$ and $\TT_{n+3}$ have Bernstein bases does not guarantee that the 2-dimension elevation procedure will be a corner cutting algorithm.

The steps to perform 1- and 2-dimension elevation can be summarized as follows.

\begin{alg}[Dimension elevation]\label{alg:order_elevation1}
Let $S(\TTS_m,\MMS,\bDelta)$ be a given piecewise Chebyshevian spline space, $\bDelta^*$ an extended partition, and $s$ a spline function of the form \eqref{eq:spline}.
To represent $s$ in the space $S(\TTS_{m+r},\widetilde{\MMS},\bDelta)$, whose dimension is elevated by $r=1,2$:
\begin{enumerate}[label=\arabic*.]
\item repeatedly apply Algorithm \ref{alg:knot_ins} to the space $S(\TTS_m,\MMS,\bDelta)$ until all the knots have multiplicity $n=m-1$;
\item on each nontrivial knot interval $[t_{\ell},t_{\ell+1}]$, determine the
  transition functions $g_i$, $i=1,\dots,m$, of $\TT_{\ell,m}$, and evaluate at
    $t_{\ell}$ the derivative of order $i$ if $r=1$, and also of order $i+1$ if $r=2$;
\item on each nontrivial knot interval $[t_{\ell},t_{\ell+1}]$, determine the
  transition functions $\tilde{g}_i$, $i=1,\dots,m+r$, of $\TT_{\ell,m+r}$, and
    evaluate at $t_{\ell}$ the derivative of order $i$
   if $r=1$, and also of order $i+1$ if $r=2$;
\item if $r=1$, determine the coefficients $\gamma_i$, $i=1,\dots,n$, through
  \eqref{eq:order_elev_1}, or, if $r=2$, determine the coefficients
    $\gamma_i$, $i=1,\dots,n$ and $\delta_i$, $i=1,\dots,n+1$, through \eqref{eq:order_elev_2a};
\item determine the coefficients $\tilde{c}_i$, $i=0,\dots,m+r$, either through \eqref{eq:coef_elev_1}
  if $r=1$ or through \eqref{eq:coef_elev_2} if $r=2$;
\item apply a knot removal algorithm to the spline space
  $S(\TTS_{m+r},\widetilde{\MMS},\bDelta)$, in which all the knots have
    multiplicity $m$, until multiplicity $\mu_i+r$ is obtained, where $\mu_i$ is the initial multiplicity of each knot in $S(\TTS_m,\MMS,\bDelta)$.
\end{enumerate}
\end{alg}

\section{Geometrically continuous Chebyshevian splines} \label{sec:gc_extension}
So far we have focused on splines where adjacent pieces are connected with parametric continuity,
which means that their derivatives at break-points must agree up to proper order of continuity.
Such splines are just an example in the wider class of \emph{geometrically continuous} splines (see \cite{DynMicchelli1989} and \cite{Barry1996,Maz2001,BCM2016} for splines based on polynomials and extended Chebyshevian spaces, respectively) and this section aims at showing that the tools developed so far also work in this more general context.

To define a space of geometrically continuous Chebyshevian splines, in addition to the quantities in Definition \ref{def:QCS},
we shall associate with each breakpoint $x_i$ a lower triangular matrix $M_i$ of size $m-\mu_i$, which we call a \emph{connection matrix}, having positive diagonal entries and first row and column equal to $(1,0,\dots,0)$. The space corresponds to Definition \ref{def:QCS} in which the parametric continuity conditions ii) are replaced by the following geometric continuity conditions:
\begin{equation}\label{eq:geo_cont}
{\rm ii}) \quad M_i \left(D^0 s_{i-1}(x_i),D^1 s_{i-1}(x_i),\dots,D^{m-\mu_i-1} s_{i-1}(x_i)\right)^T =
\left(D^0 s_i(x_i),D^1 s_i(x_i),\dots,D^{m-\mu_i-1} s_i(x_i)\right)^T,
\quad i=1,\dots,q.
\end{equation}
We denote such a space by $S(\TTS_m,\MMS,\bDelta,\bM)$, where $\bM\coloneqq(M_1,\dots,M_q)$ is the sequence of connection matrices.

Each transition functions $f_i$, $i=1,\dots,m+K$, can still be determined by Hermite interpolation, replacing
 the parametric continuity conditions \eqref{eq:cond_tf} by the geometric continuity conditions
\begin{equation}
M_j \left(D^0 f_{i,j-1}(x_{j}),D^1 f_{i,j-1}(x_{j}),\dots,D^{k_{j}} f_{i,j-1}(x_{j})\right)^T =
\left(D^0 f_{i,j}(x_{j}),D^1 f_{i,j}(x_{j}),\dots,D^{k_{j}} f_{i,j}(x_{j})\right)^T,
\; j=p_i+1,\dots,p_{i+m-1}-1.
\end{equation}
Up to this modification all the tools described in previous sections also apply to geometrically continuous splines.

\begin{rem}
In the case of knot insertion (see Section \ref{sec:knot_ins}) there are two possible scenarios. If the multiplicity of an existing breakpoint is increased, then the corresponding connection matrix must be updated by removing its last row and column. Instead, when a new breakpoint is inserted, also the corresponding connection matrix must be defined.
\end{rem}

The following example illustrates the transition functions, the B-spline basis and a curve in a space of geometrically continuous Chebyshevian splines.

\begin{figure*}[h!]
\centering
\subfigure[]{
{\includegraphics[width=0.8\textwidth/3]{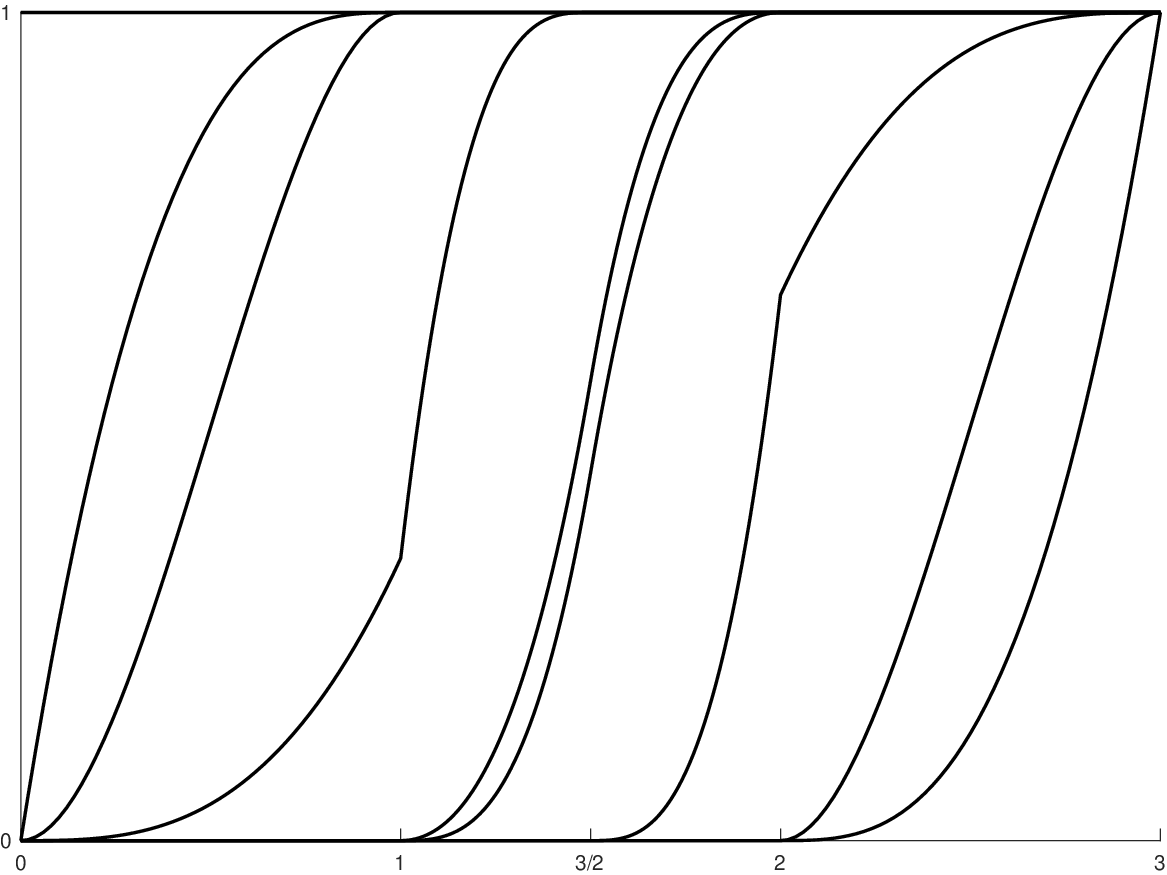}\label{fig:B-spline_geo_cont_b-7-crop}}
\hfill
}
\hfill
\subfigure[]{
{\includegraphics[width=0.8\textwidth/3]{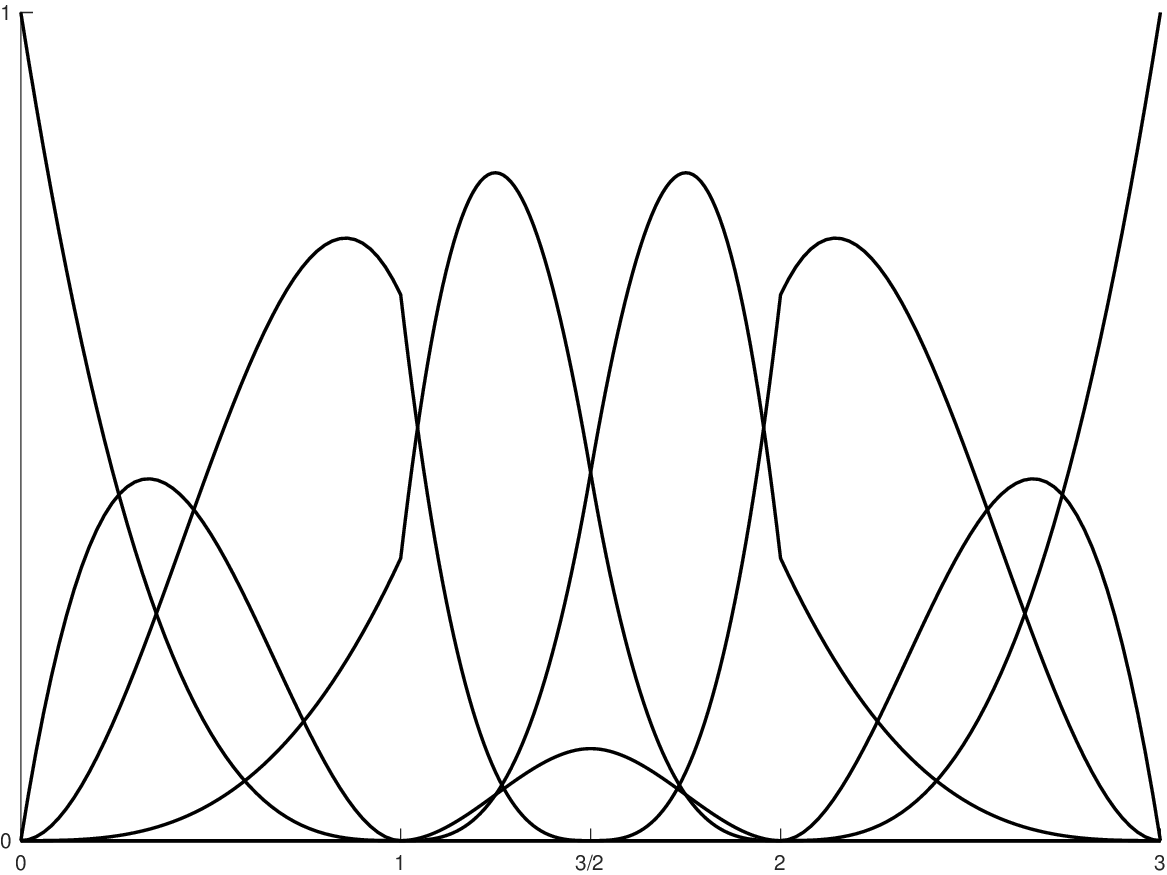}\label{fig:B-spline_geo_cont_b-7-crop}}
\hfill
}
\hfill
\subfigure[]{
{\includegraphics[width=0.8\textwidth/3]{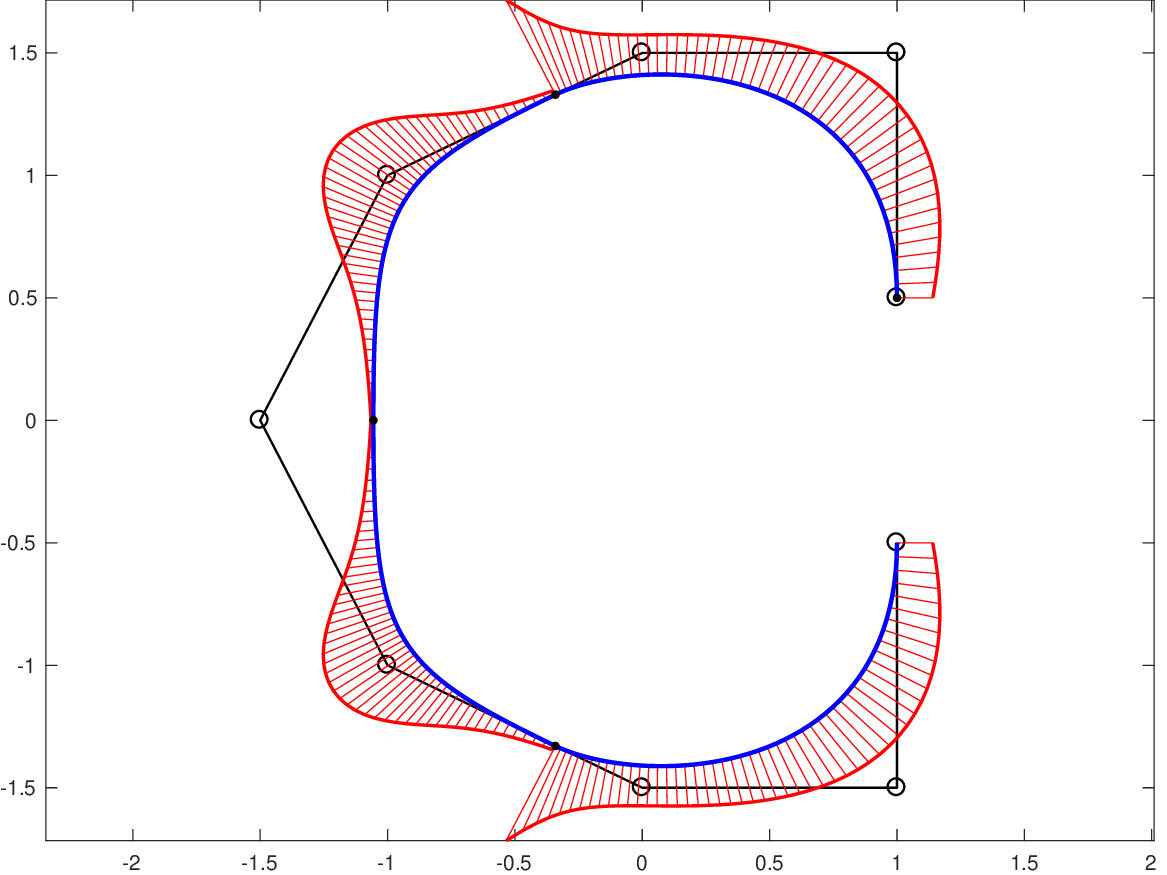}\label{fig:geo_cont_b-7-crop}}
\hfill
}
\caption{
From left to right: Geometrically continuous transition functions, corresponding B-spline basis functions and open spline curve (with curvature comb) in the spline space
of Example \ref{ex:geo_cont} with $\beta = -7$.}
\label{fig:close_geo_cont}
\end{figure*}

\begin{exmp}[Geometrically continuous Chebyshevian splines]\label{ex:geo_cont}
Let us consider the EC-spaces ${\mathcal T}_4=\Span \{1,t,\cos(t),$ $\sin(t)\}$ and $\PP_4=\Span \{1,t,\allowbreak t^2,t^3\}$, and the piecewise Chebyshevian spline space $S(\TTS_4,\MMS,\bDelta,\bM)$ defined on $[0,3]$ with $\TTS_4\coloneqq\{{\mathcal T}_4,\PP_4,\PP_4,{\mathcal T}_4\}$, break-points
$\bDelta=\{1,\frac{3}{2},2\}$, multiplicities $\MMS=(2,1,2)$ and connection matrices $\bM=(M_1,M_2,M_3)$ having the form:
\begin{equation}
M_1=
\begin{pmatrix}
1 & 0 \\
0 & 4
\end{pmatrix},
\quad
M_2=
\begin{pmatrix}
1 & 0 & 0 \\
0 & 1 & 0  \\
0 & \beta & 1
\end{pmatrix},
\quad
M_3=
\begin{pmatrix}
1 & 0 \\
0 & 1/4
\end{pmatrix},
\qquad \beta \in {\mathbb R}.
\end{equation}
Figures \ref{fig:close_geo_cont} shows the transition functions and the associated B-spline basis functions corresponding to $\beta=-7$
 and  depicts the corresponding spline curves defined by a \qtext{C}-shaped control polygon.
The multiplicities of break-points and the structure of the connection matrices make so that
the curves are symmetric, $G^1$ continuous at the break-points $x_1=1$ and $x_2=2$, and $G^2$ continuous at $x_3=\frac32$.
The superimposed curvature comb highlights the junction points with $G^1$ and $G^2$ continuities.
We emphasize that, by selecting $\beta=0$, the curve would be $C^2$ continuous at $3/2$ and $G^1$ continuous at the other knots.
It is also interesting to remark that, by selecting increasing values of the parameter
$\beta$ (which affects the curve shape in the intervals
$[1, \frac{3}{2}]$ and $[\frac{3}{2},2]$), \qtext{tension} effects can be obtained.
\end{exmp}

\section{Multi-order (multi-degree) Chebyshevian splines}\label{sec:mo-spline}
Piecewise polynomial functions whose pieces may have different degrees are called \emph{multi-degree splines}.
They were firstly introduced in \cite{SZS2003} and, more recently, B-spline bases for these spaces were constructed by means of either integral recurrence relations \cite{SW2010b,SW2010a,BCM2017}, algebraic recurrence formulas \cite{BC2019} or transition functions \cite{BCM2017}.
Alternative procedures for evaluating these splines are based on linear operators which map a known basis into the multi-degree B-spline basis \cite{TSHH2017,BC2021arxiv,BCRKIRDE} of which \cite{BC2021arxiv} is proven to be a stable algorithm. Though in the framework of polynomial splines these methods do not involve solving linear systems, they have the downside of being global, meaning that they generate all the B-spline basis functions at the same time and thus they do not lend themselves well to local modifications.

By analogy, we call \emph{multi-order (or multi-degree) Chebyshevian splines} piecewise functions defined by connecting EC-spaces of different dimensions. For these splines, an evaluation algorithm based on the aforementioned mapping between a standard polynomial basis and the multi-degree Chebyshevian B-splines has been recently proposed in \cite{chebMD2019}.
It is however important to stress that no (either theoretical or numerical) method is known at present to determine whether or not a multi-degree Chebyshevian spline space has a B-spline basis or is good for design, which inevitably limits the practical utility of the related computational methods.

In the remainder of this section we will show how to adapt the evaluation procedure based on transition functions to multi-order Chebyshevian splines. For defining such splines, denoted by $m_i$
the dimension of the EC-space on $[x_i,x_{i+1}]$, rather than the multiplicity of knots it is more natural to assign the order of continuity $k_i$ at break-points, under the condition that $0 \leq k_i < \min(m_{i-1},m_{i})$ for $i=1,\ldots,q$.
Hence, as already discussed in \cite{BMuhl2003}, for constructing the B-spline basis we shall consider two different extended partitions $\bDelta^*_s=\{s_i\}$ and $\bDelta^*_t=\{t_i\}$, defined as follows:
\[
\bDelta^*_s\coloneqq \left\{s_{1},\dots, s_{\Gamma}\right\} \equiv
        \{ \underbrace{x_0,\dots,x_0}_{m_0 \text{ times}}, \underbrace{x_1,\dots,x_1}_{m_1-k_1-1 \text{ times}}, \dots, \underbrace{x_q,\dots,x_q}_{m_q-k_q-1 \text{ times}} \},
        \]
        and
        \[
       \bDelta^*_t\coloneqq \left\{t_{1}, \dots, t_{\Gamma}\right\} \equiv
        \{ \underbrace{x_1,\dots,x_1}_{m_0-k_1-1 \text{ times}}, \dots,\underbrace{x_q,\dots,x_q}_{m_{q-1}-k_{q}-1 \text{ times}},\underbrace{x_{q+1},\dots,x_{q+1}}_{m_q \text{ times}}\},
        \]
where $\Gamma:=m_0+ \sum_{i=1}^{q} (m_i-1-k_i)$.
In this way each knot of $\bDelta_s^*$ and $\bDelta_t^*$ is respectively the starting point and endpoint of a B-spline basis function, namely the interval
$[s_i,t_i]$ is the support of the $i$th B-spline $N_i$, $i=1,\dots,\Gamma$. We also observe that
the dimension of a multi-order generalized spline space is simply the
cardinality of $\bDelta_s^*$ or of $\bDelta_t^*$.
In this setting, each transition function $f_i$, $i=2,\ldots,\Gamma$, is nontrivial in the interval $(s_i,t_{i-1})$ and can be determined from an adapted version of the conditions \eqref{eq:cond_tf}, by imposing that $f_i$ vanish $k_i^s+1$ times at $s_i$, that $1-f_i$ vanish $k_{i-1}^t+1$ times at $t_{i-1}$ and that $f_i$ have continuity $k_j$ at every breakpoint $x_j{\, \in (s_i,t_{i-1})}$, where $k_i^s$ and $k_{i-1}^t$ are defined as:
\[
k_i^s\coloneqq m_{ps_i}-\max\{j\geqslant 0\; | \; s_i=s_{i+j}\}-2\ \qquad \text{and} \qquad k_{i-1}^t\coloneqq m_{pt_i}-\max\{j\geqslant 0 \; | \; t_{i-j}=t_{i}\}-2,
\]
and where ${ps_i}$ and ${pt_i}$ are integers such ${x_{ps_i}}=s_i$ and ${x_{pt_i}} =t_i$.

\begin{rem}
Note that, also in the multi-order setting, it is possible to introduce connection matrices, thus generating multi-order Chebyshevian splines with geometric continuity similarly as for Section \ref{sec:gc_extension}, and use the transition functions for their calculation.
\end{rem}

The following example illustrates the transition functions and the corresponding B-spline basis for a parametrically continuous multi-order Chebyshevian spline space.

\begin{figure*}[h!]
\begin{center}
\hfill
\subfigure[]{\includegraphics[width=1.0\textwidth/3]{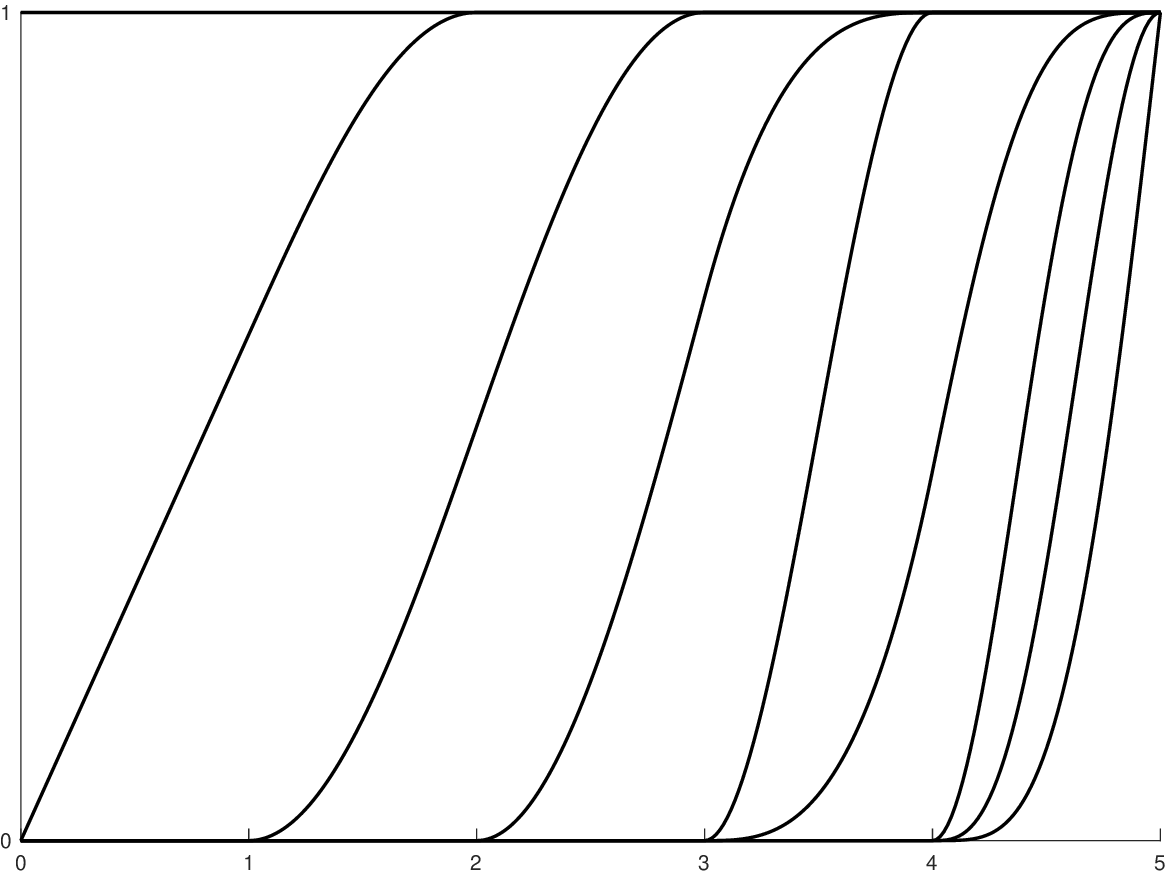}\label{fig:es4a}}\hfill
  \subfigure[]{\includegraphics[width=1.0\textwidth/3]{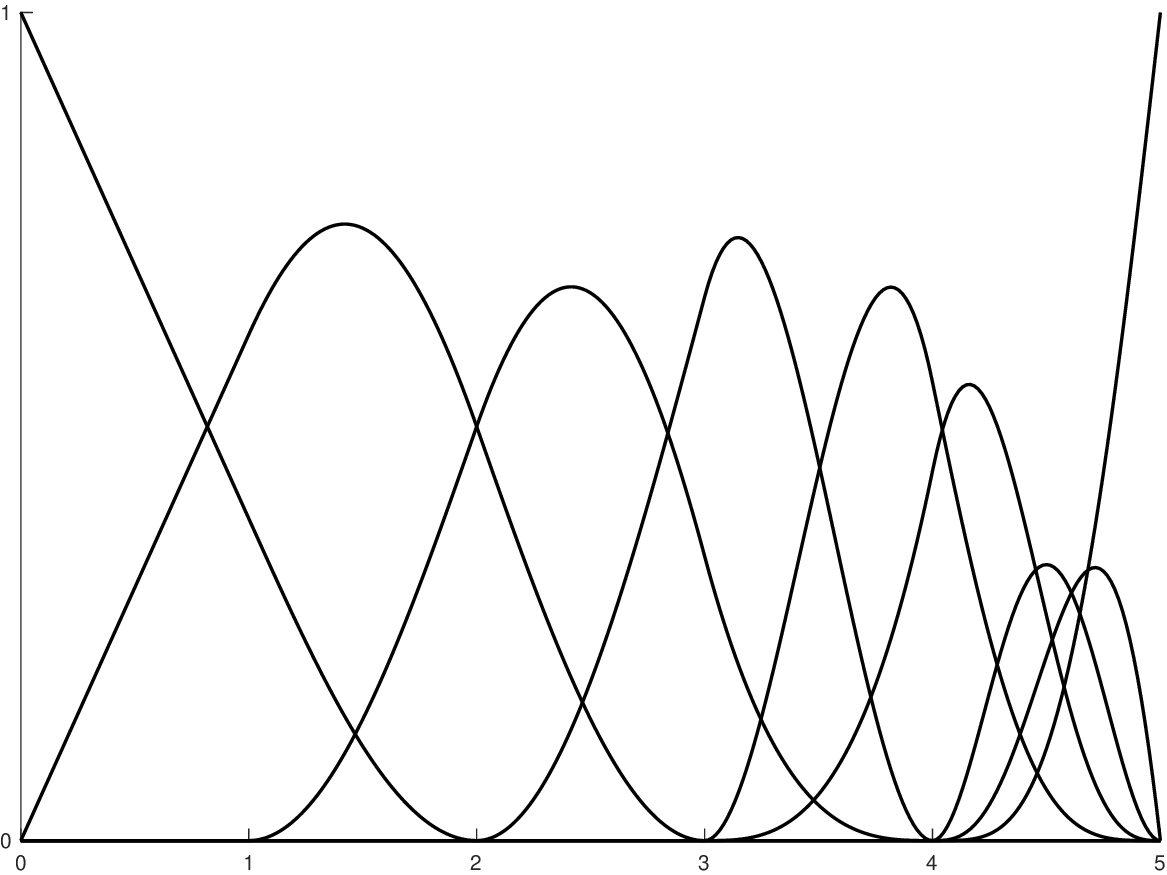}\label{fig:es4b}}\hfill
\subfigure[]{\includegraphics[width=0.8\textwidth/3,trim=3.25cm 1cm 2.75cm 0.7cm, clip=true]{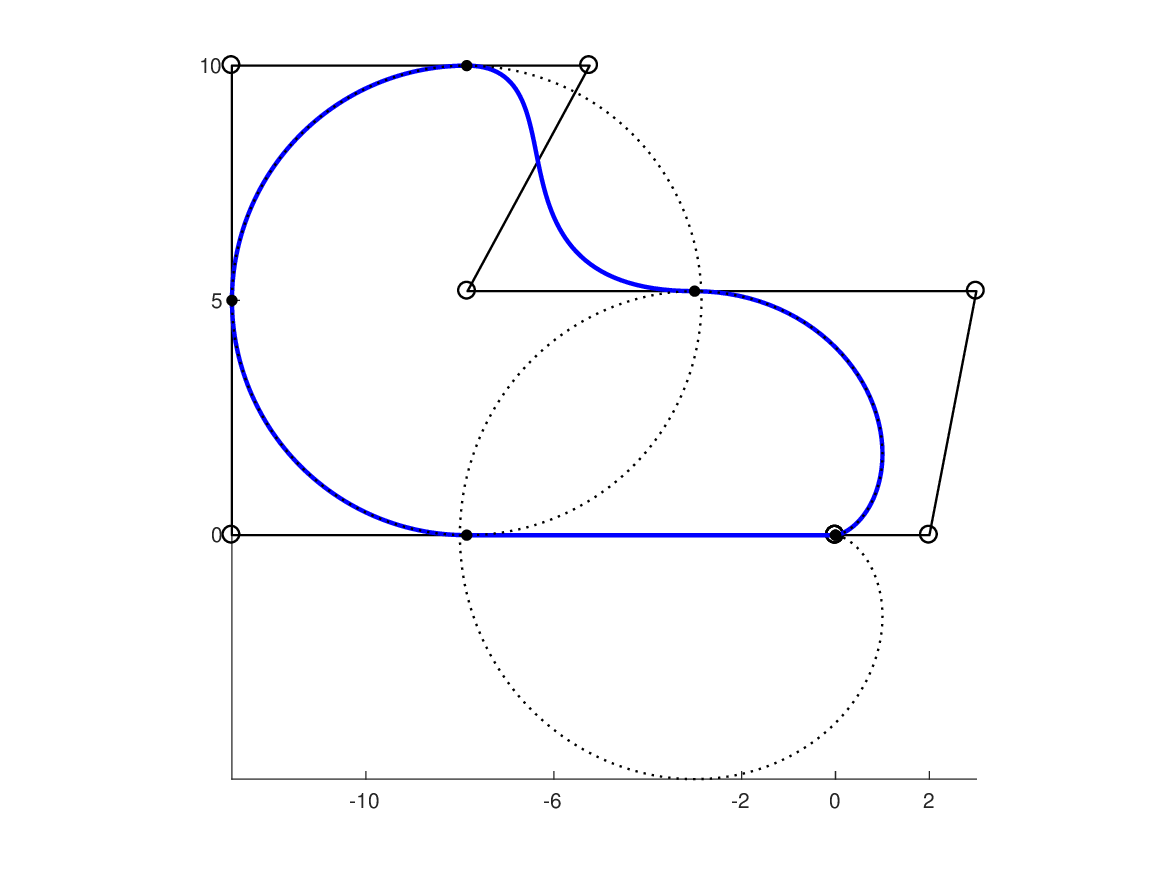}\label{fig:es4c}}
\end{center}
\vspace{-0.5cm}
  \caption{Transition functions \subref{fig:es4a},
  B-spline functions \subref{fig:es4b} and spline curve \subref{fig:es4c} in
  the multi-order spline space of Example \ref{ex:mo}.
  Note that the cardioid and the circular segments are exactly reproduced in the
  selected spline space.}
\label{fig:close_mo}
\end{figure*}

\begin{exmp}[Multi-order Chebyshevian spline]\label{ex:mo}
Multi-order Chebyshevian splines can be used to represent special curves by assigning on each interval the EC-space of the minimum possible dimension.
To illustrate this fact, Figure \ref{fig:close_mo} shows a closed, nonperiodic curve consisting of five pieces, which include a straight line segment, two circular arcs, a cubic polynomial and a cardioid segment.
Such a curve has been obtained by taking a spline space with $[a,b]=[0,5]$, clamped uniform knot partition, $C^1$ continuity at breakpoints and the following sequence of section spaces
$$
\begin{array}{c}
\TT_{0,2}=\Span \{1,t \}, \quad  \TT_{1,3}=\TT_{2,3}=\Span \left \{1, \cos\left(\frac{\pi}{2} t \right), \sin\left(\frac{\pi}{2} t \right) \right \}, \quad
\TT_{3,4}=\Span \{1,t,t^2,t^3\}, \\
\TT_{4,5}=\Span \{1, \cos(\phi t), \sin(\phi t), \cos(2\phi t), \sin(2\phi t)\} \quad \hbox{with} \quad \phi=\frac{2}{3}\pi.
\end{array}
$$
\end{exmp}

\section{Comparisons with the evaluation method relying on the extraction operator}\label{sec:comparison}

It is certainly well-known that stability problems might arise when evaluating the Bernstein basis of an EC-space, also known as Chebyshevian Bernstein basis \cite{AGO2019}. As already discussed in this recent paper, such
stability problems are very frequently encountered when considering, for instance:
\begin{itemize}
\item large shape parameters $\phi$ in EC-spaces $\TT_{i,m}$ whose set of generators contains the hyperbolic functions $\cosh(\phi x)$, $\sinh(\phi x)$, so that $\phi(x_{i+1}-x_i)\gg 0$;
\item high-dimensional EC-spaces, so that $m \gg 0$ with $m={\rm dim}(\TT_{i,m})$;
\item large-width intervals $[x_i, x_{i+1}]$, so that $x_{i+1}-x_i \gg 0$.
\end{itemize}

Since sequences of EC-spaces are building blocks of Chebyshevian spline spaces, it is clear that similar or even
worse problems of numerical instability might be encountered when evaluating Chebyshevian splines. This has been
recently observed in a manuscript that was registered in the arXiv repository a few months ago \cite{SPE2020}, with the purpose
of introducing a practical framework to deal with Chebyshevian splines. In particular, in this manuscript a
computational procedure to represent all parametrically continuous Chebyshevian splines as a linear combination of
Chebyshevian Bernstein functions related to each of the local EC-spaces in the considered sequence, has been described.
This procedure relies on the so-called \emph{extraction operator} - an operator that maps the B-spline basis of a Chebyshevian
spline space to a local Chebyshevian Bernstein basis \cite{chebMD2019} - and can be summarized as follows.\\

{\sl Evaluation Algorithm relying on the Extraction Operator (for short, EOeval):}
\begin{itemize}
\item EOeval.1) for each interval $[x_i, x_{i+1}]$, evaluate $u_{i,1},\dots,u_{i,m}$, the generators of the EC-space $\TT_{i,m}$, to determine the Wronski matrices
at the endpoints $x_i$ and $x_{i+1}$;
\item EOeval.2) for each interval $[x_i, x_{i+1}]$, compute the Bernstein basis $B_{0,m-1},\dots,B_{m-1,m-1}$ of the EC-space $\TT_{i,m}$ by solving as many linear
systems as the EC-space dimension, i.e. $m$;
\item EOeval.3) for each interval $[x_i, x_{i+1}]$, evaluate the computed Bernstein basis functions to determine the Wronski matrices at
the endpoints $x_i$ and $x_{i+1}$, with the final purpose of constructing the matrix that stores the continuity conditions
of the Chebyshevian spline space;
\item EOeval.4) compute the Extraction Operator according to Algorithm \cite[Fig 1]{SPE2020}, which entails determining (for each internal break-point and for each continuity
order) the null space of a suitable vector according to Algorithm \cite[Fig 2]{SPE2020}.
\end{itemize}

In this section we aim at comparing the evaluation algorithm relying on the Extraction Operator (for short, EOeval) with our evaluation
algorithm relying on the transition functions (shortly denoted with the acronym TFeval).
Thus, for the sake of clarity, we also provide the pseudo code of the algorithm we propose
for the evaluation of the B-spline basis of a Chebyshevian spline space.\\

{\sl Evaluation Algorithm relying on the Transition Functions (for short, TFeval):}
\begin{itemize}
\item TFeval.1) for each interval $[x_i, x_{i+1}]$, evaluate the generators $u_{i,1},\dots,u_{i,m}$ of the EC-space $\TT_{i,m}$ to determine the Wronski matrices at the endpoints $x_i$ and $x_{i+1}$;
\item TFeval.2) solve as many linear systems as ${\rm dim}(S)-1$ (where $S$ denotes the Chebyshevian spline space), to
determine all the needed transition functions $f_i$, $i = 2, \ldots, {\rm dim}(S)-1$.
\end{itemize}

Observing the two pseudo-codes, it is evident that the first step is exactly the same.
Thus, in order to have a fair comparison, in the first step of both the evaluation algorithms the same expressions for the
generators of the EC-space $\TT_{i,m}$ will be used and, among all possible equivalent choices, {the most stable one} will be selected (see again \cite{AGO2019,SPE2020}).

As to the second step, each algorithm requires the solution of several linear systems. More precisely, with respect to EOeval, the TFeval algorithm has to solve fewer but higher-dimensional linear systems. Just to give a concrete idea, for a Chebyshevian spline space of dimension $10$ with $m=7$ and $K=q=3$, we have that the TFeval algorithm has to solve just $9$ linear systems having size between $7 \times 7$ and $28 \times 28$, whereas the EOeval algorithm has to solve $28$ linear systems all having size $7 \times 7$.
But, while after solving these linear systems TFeval has to do nothing else, Eoeval has to perform two more steps to reach the sought result.
It is also important to underline that, even if the linear systems to be solved by the two algorithms differ a lot in number and size, there are no substantial differences as far as conditioning is concerned.
However, all the examples included in the following show that, although the condition numbers of the linear systems solved by the two methods are approximately of the same order, the linear systems solved by TFeval are less sensible to rounding errors on input data and lead to  basis functions that are more accurate as well as capable of preserving the symmetry properties of the spline space.\\
In all the considered examples, in order to compare the two evaluation algorithms from the point of view of numerical stability and accuracy, we have to calculate the numerical errors provided by each method.
Thus, we assume that the “exact” reference values are the ones obtained by evaluating the B-spline basis starting from a symbolic representation matrix provided by MATLAB's Symbolic Math Toolbox, whereas the numerical results rely on  MATLAB's variable-precision floating-point arithmetic (VPA) with 32-digit precision.\\
The graphs accompanying the examples show the maximum absolute error between the evaluations of all the B-spline basis functions (calculated in 32-digit precision on a discrete number of points) and the exact reference values.
For the sake of brevity, during the description of the examples, we will say that these error values are obtained by applying the \emph{Symbolic Error Test} to the two evaluation methods.\\
In addition, in the case of symmetric Chebyshevian spline spaces, we also include the maximum absolute error obtained by
comparing (on a certain number of evaluation points) the two basis functions that are known to be symmetric about the center of the interval.
As to the error values obtained with this second test, we will say that they are obtained by applying the \emph{Symmetry Check} to the two evaluation methods.\\
In order to consider a selection of spaces with very different features, the selected examples range from the case of a Chebyshevian spline space that reduces to a single EC-space (Examples \ref{ex:EC1} and \ref{ex:EC2}), to the more representative cases of mono-order (Example \ref{ex:monoorder}) and multi-order (Example \ref{ex:multiorder}) Chebyshevian spline spaces.

\begin{exmp}\label{ex:EC1}
Let $[0,4]$ be the domain of definition, and consider the EC-spaces
\begin{equation} \label{eq:iperb}
  \Span\{1,\ldots, x^n,\cosh(10 x),\sinh(10 x)\}  \quad \hbox{with} \quad n=1,\ldots,14.
\end{equation}
To compare the Chebyshevian Bernstein bases computed by the two evaluation algorithms (TFeval and EOeval)
for each of the above EC-spaces, we apply first the Symbolic Error Test and then, for each value of $n$, we display the largest error provided by each method (see Figure \ref{fig:es1a}).

As for the EC-space of dimension $16$ (corresponding to the choice $n=13$), we show in Table \ref{tab:es1_symbolic} all the partial results of the Symbolic Error Test (i.e., the absolute errors regarding each basis function). Note that the maximum of the error values shown in this table is exactly the largest error shown for $n=13$ in the graph of Figure \ref{fig:es1a}.

For the same EC-space, we also show in Table \ref{tab:es1_symmetry} the result of the Symmetry Check, i.e., the maximum absolute errors
obtained by comparing the two basis functions that are known to be symmetric about the center of the interval.
Furthermore, in Figure \ref{fig:es1_basi} we illustrate the graphs of the Bernstein basis functions provided by the two evaluation algorithms (TFeval and EOeval). As it is evident, the second and the second to last basis functions (see Table \ref{tab:es1_symmetry}, $i=1,14$) computed by the EOeval algorithm are visually very
different (as it is natural to expect if the absolute error between the two functions is of order $10^{-1}$), thus confirming the lower accuracy of EOeval with respect to TFeval.

\begin{table}[h!]
\begin{center}
\begin{tabular}{rcc}
\hline
  $B_{i,15}$ &   \textrm{TFeval} & \textrm{EOeval} \\
\hline
  0 & 6.661338147750939e-16 & 9.630723062092894e-07 \\
  1 & 3.894329303477662e-12 & 2.607500630625010e-01 \\
  2 & 2.635047735566332e-11 & 1.826763789509578e-03 \\
  3 & 1.680284800187337e-10 & 1.067719375855147e-03 \\
  4 & 3.093084677274760e-10 & 5.952836528486105e-05 \\
  5 & 3.497080403036534e-10 & 7.395153359374129e-06 \\
  6 & 2.911615393230704e-10 & 7.177414484693667e-07 \\
  7 & 1.493716261791178e-10 & 8.073124388441322e-08 \\
  8 & 4.632472183629943e-11 & 4.784319329598219e-09 \\
  9 & 6.892264536872972e-12 & 5.245952283683408e-10 \\
  10 & 1.891375944751417e-12 & 3.099415168961173e-11 \\
  11 & 5.378475442796571e-13 & 4.545669396449625e-12 \\
  12 & 1.270442084866374e-13 & 1.513789094076401e-13 \\
  13 & 6.847300504375653e-14 & 6.106226635438361e-15 \\
  14 & 3.219646771412954e-15 & 9.992007221626409e-16 \\
  15 & 5.551115123125783e-17 & 4.135903062765138e-25 \\
\hline
\end{tabular}
\caption{Absolute errors provided by the Symbolic Error Test for all Chebyshevian Bernstein functions of the EC-space in \eqref{eq:iperb} with $n=13$.}
\label{tab:es1_symbolic}
\end{center}
\end{table}

\begin{table}[h!]
\begin{center}
\begin{tabular}{lcc}
\hline
  $B_{i,15}$ &   \textrm{TFeval} & \textrm{EOeval} \\
\hline
  0,15&    2.109423746787797e-15&    9.626665849749028e-07 \\
  1,14&    3.896674649617182e-12&    2.607414289674860e-01 \\
  2,13&    2.640754281912905e-11&    1.822789713677941e-03 \\
  3,12&    1.681335071168633e-10&    1.067284258902568e-03 \\
  4,11&    3.097465062218419e-10&    5.949644951949717e-05 \\
  5,10&    3.498862866102570e-10&    7.393873874966816e-06 \\
  6,9&    2.949248900652179e-10&    7.171351926738012e-07 \\
  7,8&    1.825976592151335e-10&    8.544727658765794e-08 \\
\hline
\end{tabular}
\caption{Absolute errors provided by the Symmetry Check for all Chebyshevian Bernstein functions of the EC-space in \eqref{eq:iperb} with $n=13$.}
  \label{tab:es1_symmetry}
\end{center}
\end{table}

\begin{exmp}\label{ex:EC2}
We consider the EC-space \eqref{eq:iperb} with $n=5$ (i.e., dimension 8) and, for domains $[0,b]$ of increasing width (specifically for $b \in \{ 10^{-3}, 10^{-2}, 10^{-1}, 1, 2, 2^2, 2^3, 2^4, 2^5, 2^6 \}$), we apply the Symbolic Error Test and show in Figure \ref{fig:es1b} the maximum absolute error achieved by each evaluation method. As it appears, the TFeval algorithm is again more accurate than its competitor.
\end{exmp}

\begin{figure}[h!]
\begin{center}
\hfill
\subfigure[]{\includegraphics[width=0.45\textwidth]{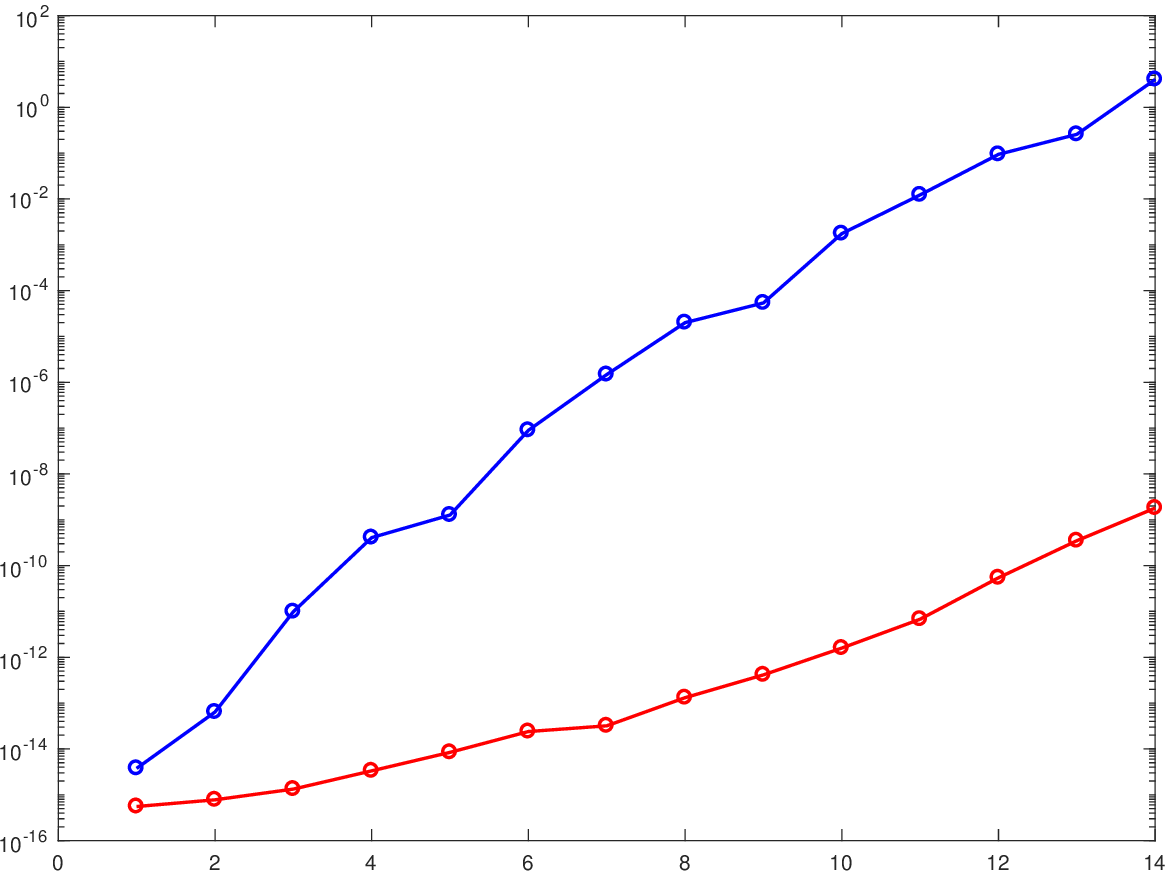}\label{fig:es1a}}\hfill
\subfigure[]{\includegraphics[width=0.45\textwidth]{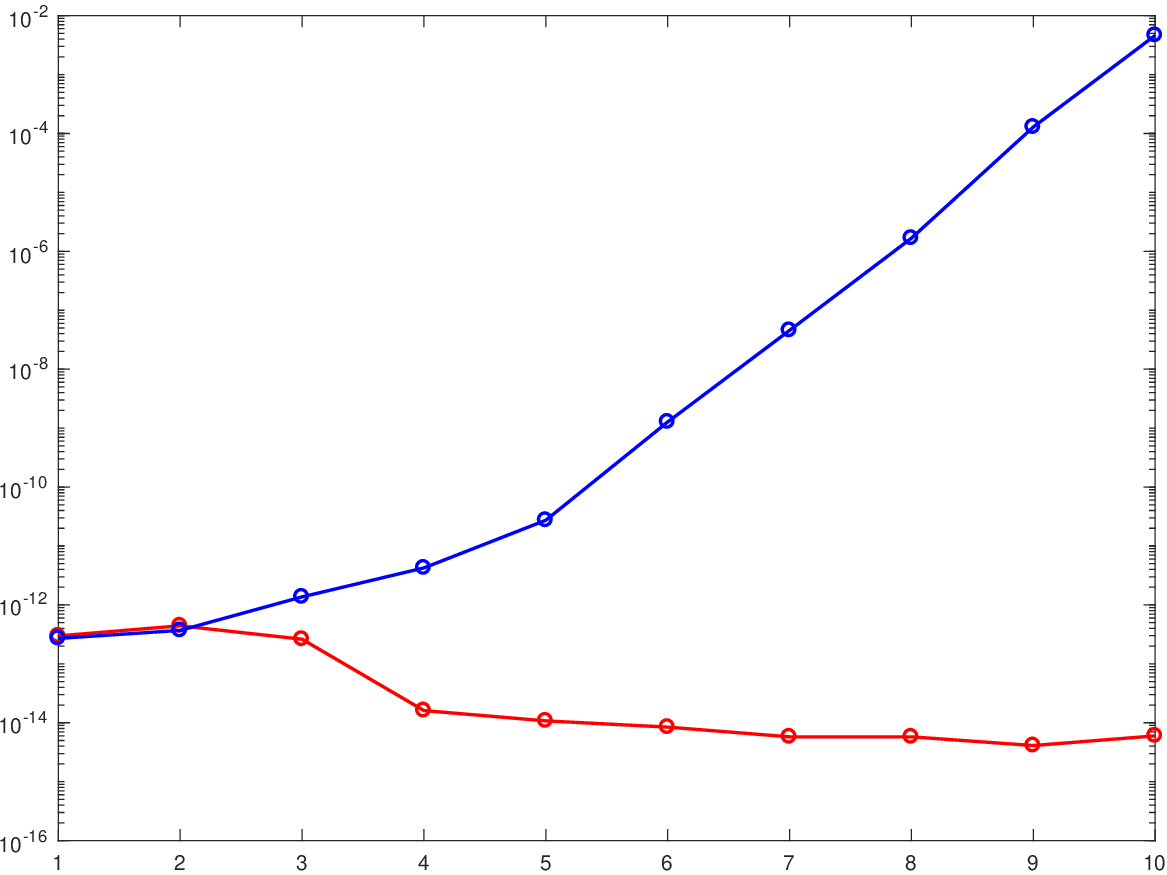}\label{fig:es1b}}\hfill
\end{center}
\vspace{-0.5cm}
  \caption{Maximum absolute errors provided by TFeval (red) and EOeval (blue) when the Symbolic Error Test is applied to the Bernstein basis of: \subref{fig:es1a} an EC-space of increasing dimension (see Example \ref{ex:EC1}) and \subref{fig:es1b}  an EC-space defined on intervals of increasing width (see Example \ref{ex:EC2}).}
\label{fig:es1}
\end{figure}

\begin{figure}[h!]
\begin{center}
\hfill
\subfigure[]{\includegraphics[width=0.45\textwidth]{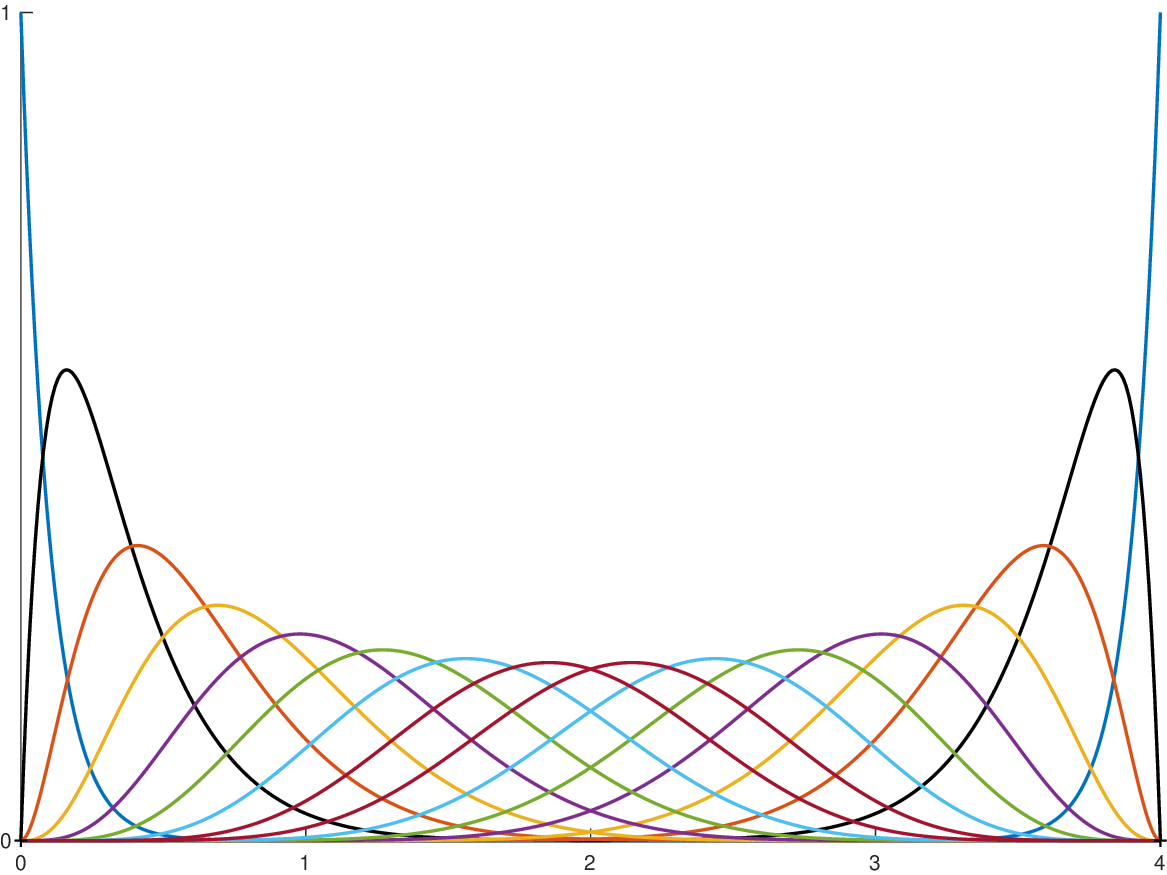}\label{fig:es1a_tf}}\hfill
\subfigure[]{\includegraphics[width=0.45\textwidth]{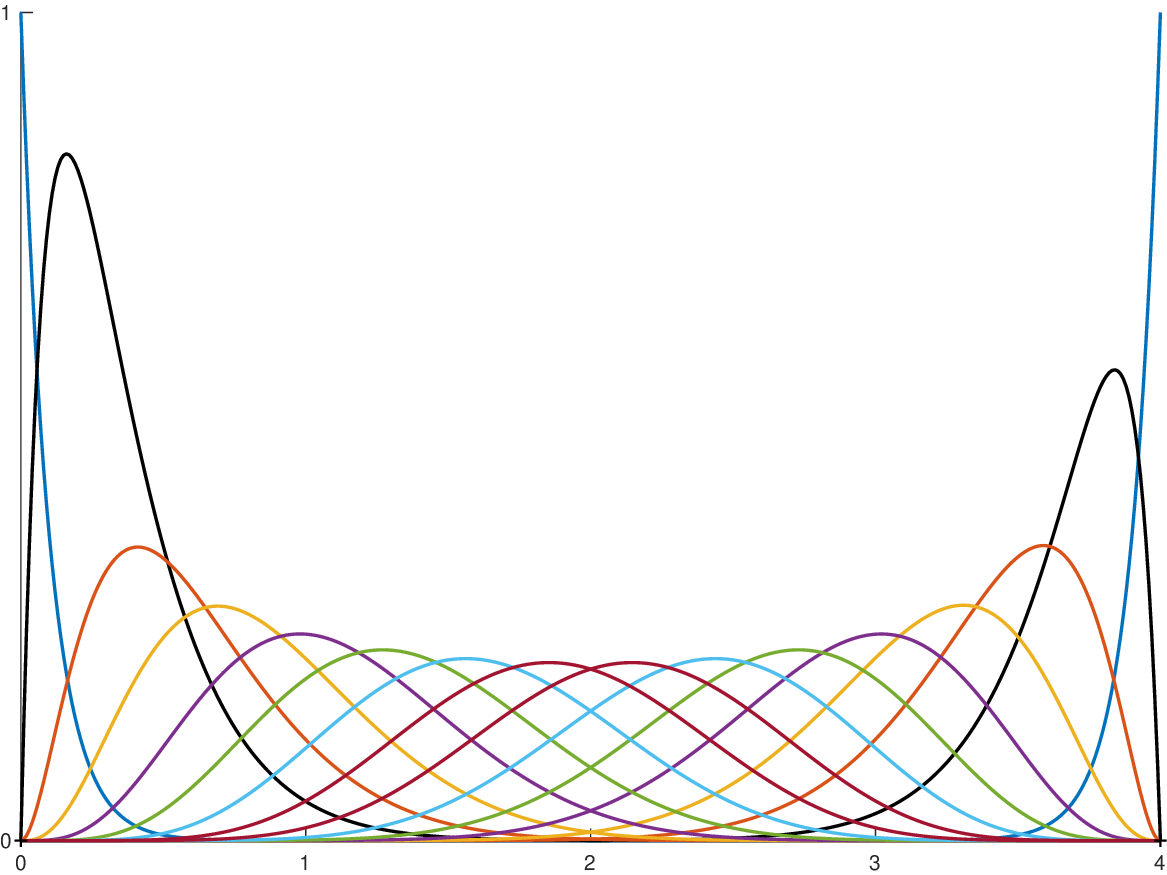}\label{fig:es1a_eo}}\hfill
\end{center}
\vspace{-0.5cm}
  \caption{Chebyshevian Bernstein basis of  the EC-space in \eqref{eq:iperb} with $n=13$: basis functions provided by TFeval \subref{fig:es1a_tf} and by EOeval \subref{fig:es1a_eo}.}
\label{fig:es1_basi}
\end{figure}
\end{exmp}

\begin{figure}[h!]
\begin{center}
\hfill
\subfigure[]{\includegraphics[width=0.45\textwidth]{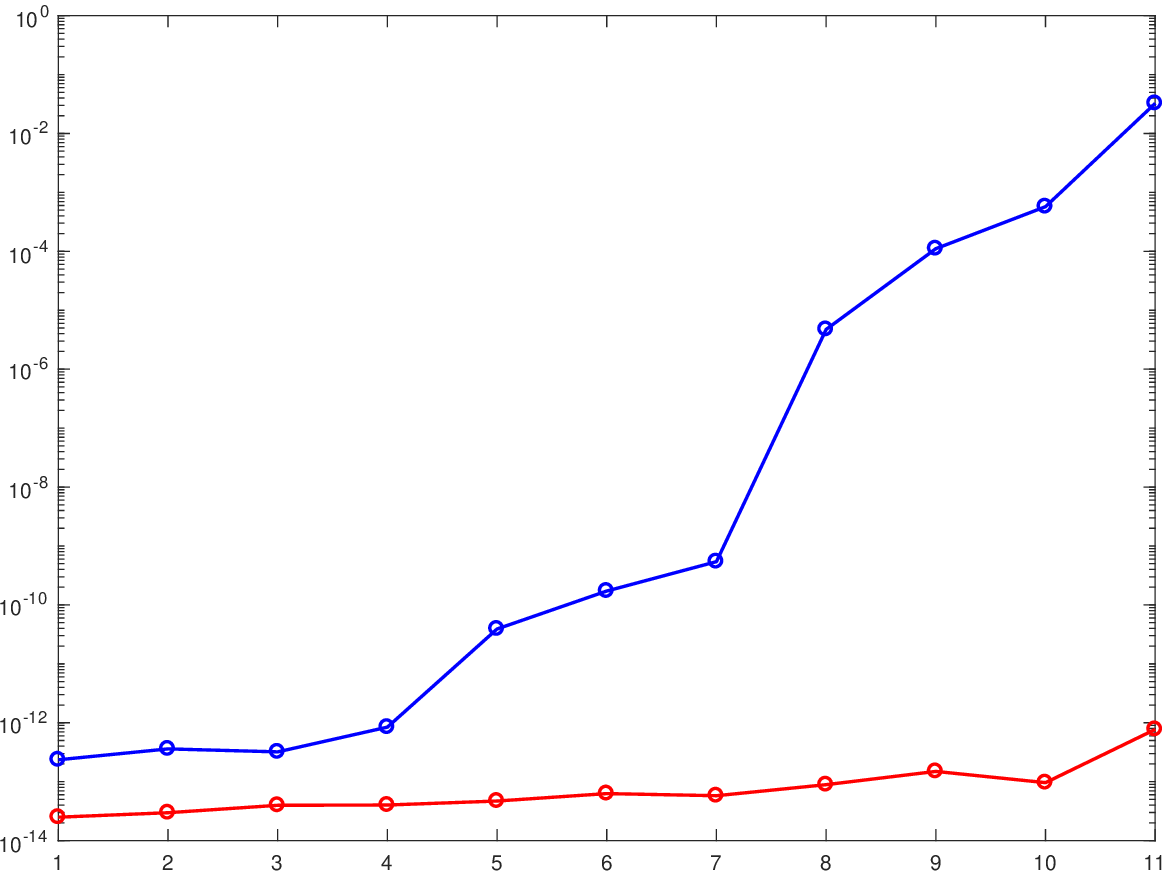}\label{fig:es2a}}\hfill
\subfigure[]{\includegraphics[width=0.45\textwidth]{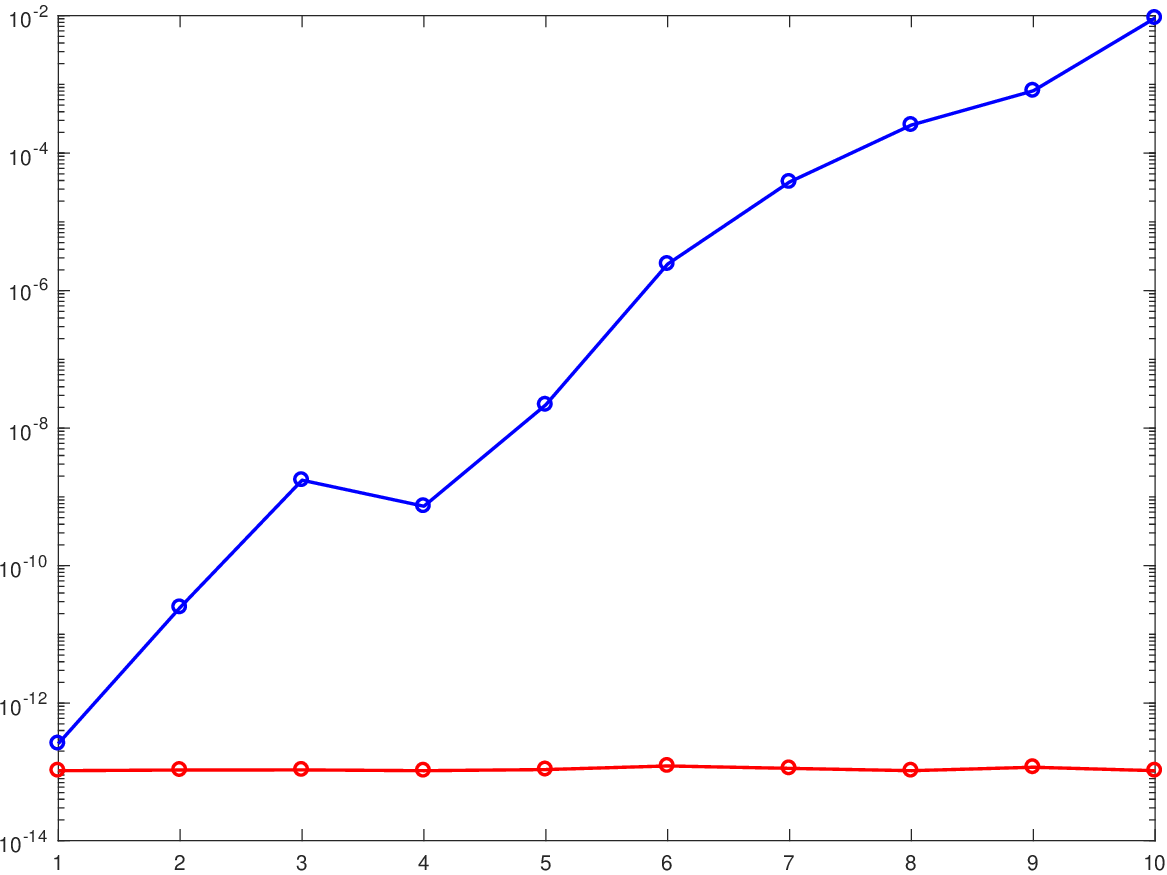}\label{fig:es2b}}\hfill
\end{center}
\vspace{-0.5cm}
\caption{
Maximum absolute errors provided by TFeval (red) and EOeval (blue) when the Symbolic Error Test is applied to the B-spline basis of: \subref{fig:es2a} the Chebyshevian spline space of Example \ref{ex:monoorder} and \subref{fig:es2b} the Chebyshevian spline space of Example \ref{ex:multiorder}.}
\label{fig:es2}
\end{figure}

\begin{exmp}\label{ex:monoorder}
We consider the sequence $\TTS_8=\{\TT_{0,8},\TT_{1,8},\TT_{2,8},\TT_{3,8}\}$ consisting of the four 8-dimensional EC-spaces of the form
$$
  \TT_{0,8} = \TT_{3,8} = \Span\{1,\ldots, x^5,\cos(x),\sin(x)\}, \qquad
  \TT_{1,8} = \TT_{2,8} = \Span\{1,\ldots, x^5,\cosh(x),\sinh(x)\}.
$$
Then, for $w \in \{0.5,0.6,0.7,0.8,0.9,0.95,0.96,0.99,0.995,0.996,0.999\}$, we construct a
Chebyshevian spline space $S(\TTS_8,\MMS,\bDelta)$ defined on $[0,2]$ with
break-points $\bDelta=\{1-w,1,1+w\}$ and multiplicities $\MMS=(1,1,1)$ (i.e., global continuity $C^6$).
Clearly, by modifying the value of $w$ in the considered set, the break-points in the partition $\bDelta$ become more and more unevenly spaced.
To compare the Chebyshevian B-spline bases computed by the two evaluation algorithms (TFeval and EOeval)
for each of the above spline spaces, we apply first the Symbolic Error Test and then, for each value of $w$, we display the largest error provided by each method (see Figure \ref{fig:es2a}).\\
Since the considered spline spaces are all symmetric, for the largest value of $w$ we also show in Table \ref{tab:es2} the result of the Symmetry Check, i.e., the maximum absolute errors
obtained by comparing the two basis functions that are known to be symmetric about the center of the interval.

\begin{table}[h!]
\begin{center}
\begin{tabular}{lcc}
\hline
  $N_{i,11}$ &  \textrm{TFeval} & \textrm{EOeval} \\
\hline
  1,11&    2.738365090237949e-13&    1.503575042249850e-12 \\
  2,10&    2.733369086627135e-13&    1.917160874498336e-12 \\
  3,9&    2.201017146319373e-14&    9.117262500524248e-12 \\
  4,8&    5.154904281212680e-14&    3.239101038779496e-02 \\
  5,7&    6.734387793781455e-14&    3.124764433409294e-02 \\
  6,6&    3.025357742103552e-14&    2.540930784666884e-06 \\
\hline
\end{tabular}
  \caption{Maximum absolute errors provided by applying the Symmetry Check to the Chebyshevian spline space $S(\TTS_8,\MMS,\bDelta)$ of Example \ref{ex:monoorder} with $\bDelta$ obtained from the largest value of $w$.}
  \label{tab:es2}
\end{center}
\end{table}

\end{exmp}

\begin{exmp}\label{ex:multiorder}
We consider the sequence $\{\TT_{0,8},\TT_{1,6},\TT_{2,8},\TT_{3,6},\TT_{4,8}\}$ consisting of the EC-spaces of the form
$$
\begin{array}{cc}
  \TT_{0,8} = \TT_{4,8} = \Span\{1,\ldots, x^5,\cos(x),\sin(x)\}, \quad
  \TT_{2,8} = \Span \{1,\ldots, x^7\}, \\
  \TT_{1,6} = \TT_{3,6} = \Span \left\{1,\ldots, x^3,\cosh\left(\frac{3}{2} x \right),\sinh\left(\frac{3}{2} x \right) \right\}.
\end{array}
$$
Then, for $w \in \{0.5,0.9,0.95,0.96,0.97,0.99,0.995,0.996,0.997,0.998\}$, we construct the multi-order Chebyshevian spline space defined on
$[0,5]$ by the break-points $\bDelta=\left\{\frac{1}{2}, 3-\frac{5}{2}w, 2+\frac{5}{2}w, \frac{9}{2} \right\}$ and by the global order of continuity $C^5$ (i.e., $k_i=5 \; \forall i$).
Clearly, by modifying the value of $w$ in the considered set, the break-points in the partition $\bDelta$ become more and more unevenly spaced.
To compare the multi-order Chebyshevian B-spline bases computed by the two evaluation algorithms (TFeval and EOeval)
for each of the above spline spaces, we apply first the Symbolic Error Test and then, for each value of $w$, we display the largest error provided by each method (see Figure \ref{fig:es2b}).
\end{exmp}

The error trend in Figures \ref{fig:es1}, \ref{fig:es2} also shows that TFeval can successfully
handle Chebyshevian spline spaces with much higher dimensions and much more uneven knot partitions of those considered in the above experiments, to the contrary of EOeval.

\section{Conclusion}

We have presented a method for computing the B-spline basis of a piecewise Chebyshevian spline space good for design, which gives rise to an efficient and accurate evaluation algorithm. The proposed approach is based on the computation of a set of functions, called the transition functions, which form a basis for the space of interest and each of which can be calculated independently of the others by solving a small-size linear system. As a consequence of its local nature, the method lends itself well to handle local modifications of the data, as well as to efficiently perform knot insertion and order (or degree) elevation. We have shown that this computational approach is more stable and accurate than the only other available algorithm to date, both in the classical setting — in which all section spaces have the same dimension — and in the case where the section spaces may have different dimensions.
The comparison of the two methods also emphasizes that our proposal can successfully handle spline spaces with much higher dimensions and much more uneven knot partitions.


\section*{Acknowledgements}
The authors are grateful to Marie-Laurence Mazure for her comments and useful suggestions. All the authors are members of the Italian GNCS-INdAM which has partially supported this work. The authors are also members of the Alma Mater research center on Applied Mathematics ($AM^2$).




\ifelsevier
\bibliographystyle{model3-num-names} 
\bibliography{spline_trans_fun_partII} 
\fi



\ifspringer
\bibliographystyle{spmpsci} 
\bibliography{spline_trans_fun_partII} 
\fi

\typeout{get arXiv to do 4 passes: Label(s) may have changed. Rerun}

\end{document}